\documentclass[12pt]{sn-jnl}% Math and Physical Sciences Reference Style
\theoremstyle{thmstyleone}%
\theoremstyle{thmstyletwo}%

\theoremstyle{thmstylethree}%
\renewcommand{\Re}{\mathbb{R}}

\newcommand{\wv}{\mathbf{w}}
\newcommand{\av}{\mathbf{a}}
\newcommand{\bv}{\mathbf{b}}
\newcommand{\cv}{\mathbf{c}}
\newcommand{\dv}{\mathbf{d}}

\newcommand{\sv}{\mathbf{s}}

\newcommand{\uv}{\mathbf{u}}
\newcommand{\vv}{\mathbf{v}}

\newcommand{\xv}{\mathbf{x}}
\newcommand{\yv}{\mathbf{y}}
\newcommand{\zv}{\mathbf{z}}
\newcommand{\At}{\mathcal{A}}

\newcommand{\Pt}{\mathcal{P}}

\newcommand{\Xt}{\mathcal{X}}
\newcommand{\Yt}{\mathcal{Y}}

\newcommand{\Rm}{\mathbf{R}}
\newcommand{\one}{\mathbf{e}}
\newcommand{\prob}{\mathrm{Prob}}

\newtheorem{Thm}{Theorem}[section]
\newtheorem{Def}[Thm]{Definition}

\newtheorem{Lem}[Thm]{Lemma}
\newtheorem{Cor}[Thm]{Corollary}
\newtheorem{Exm}[Thm]{Example}

\title[Multi-Linear Pseudo-PageRank for Hypergraph Partitioning]{Multi-Linear Pseudo-PageRank for Hypergraph Partitioning}

\author[1]{\fnm{Yannan} \sur{Chen}}\email{ynchen@scnu.edu.cn}

\author[1]{\fnm{Wen} \sur{Li}}\email{liwen@scnu.edu.cn}
\author*[2]{\fnm{Jingya} \sur{Chang}}\email{jychang@gdut.edu.cn}
\affil*[1]{\orgdiv{School of Mathematical Sciences}, \orgname{South China Normal University}, \orgaddress{\city{Guangzhou}, \postcode{510631}, \country{China}}}

\affil[2]{\orgdiv{School of Mathematics and Statistics}, \orgname{Guangdong University of Technology}, \orgaddress{\city{Guangzhou}, \postcode{510520}, \country{China}}}

\abstract{Motivated by the PageRank model for graph partitioning, we develop an extension of PageRank for partitioning uniform hypergraphs. Starting from adjacency tensors of uniform hypergraphs, we establish the multi-linear pseudo-PageRank (MLPPR) model, which is formulated as a multi-linear system with nonnegative constraints. The coefficient tensor of MLPPR is a kind of Laplacian tensors of uniform hypergraphs, which are almost as sparse as adjacency tensors since no dangling corrections are incorporated. Furthermore, all frontal slices of the coefficient tensor of MLPPR are M-matrices. Theoretically, MLPPR has a solution, which is unique under mild conditions. An error bound of the MLPPR solution is analyzed when the Laplacian tensor is slightly perturbed. Computationally, by exploiting the structural Laplacian tensor, we propose a tensor splitting algorithm, which converges linearly to a solution of MLPPR. Finally, numerical experiments illustrate that MLPPR is powerful and effective for hypergraph partitioning problems.}


\begin{document}
\maketitle

%\begin{equation*}\label{prob}
%  \alpha
%\end{equation*}

%\end{document}

\section{Introduction}

PageRank is the first and best-known method developed by Google to evaluate the importance of web-pages via their inherent graph and network structures \cite{BP'98, BL'06}. There are three advantages of PageRank: existence and uniqueness of the PageRank solution, fast algorithms for finding the PageRank solution, and a built-in regularization \cite{Gl'15}. Thus, PageRank and its generalizations have been extensively used in network analysis, image sciences, sports, biology, chemistry, mathematical systems, and so on \cite{Gl'15, LM'06}. Hypergraphs have the capability of modeling connections among objects according to their inherent multiwise similarity and affinity. Recently, higher-order organization of complex networks at the level of small network subgraphs becomes an active research area \cite{BGL16, BGL'17, KW06, RELWL14}. The set of small network subgraphs with common connectivity patterns in a complex network could be modeled by a uniform hypergraph. A natural question is how can we generalize the PageRank idea from graphs to uniform hypergraphs.

PageRank was originally interpreted as a special Markov chain \cite{LM'06}, where the transition probability matrix is a columnwise-stochastic matrix. Many algorithms, e.g., power methods, were employed to find the dominant eigenvector of the PageRank eigensystem \cite{IS08}. Alternatively, PageRank could be represented as a linear system with an M-matrix, which has favorable theoretical character \cite{F15} and good numerical performance \cite{LM06}. In the context of PageRank for partitioning a graph, although a Laplacian matrix of the graph is an M-matrix, an adjacency matrix of the graph is not necessarily columnwise-stochastic, when the graph has dangling vertices \cite{Ch'05}. Some traditional approaches \cite{ACL'07,ACL'08,Gl'15,LM06} incorporated a rank-one correction into a graph adjacency matrix and thus produced a columnwise-stochastic matrix. Using this columnwise-stochastic matrix, PageRank generates the stochastic solution. However, the rank-one correction of the graph adjacency matrix has no graphic meaning. Alternatively, following the framework of M-linear systems, Gleich \cite{Gl'15} used the graph Laplacian matrix, which is only columnwise-substochastic, to construct a pseudo-PageRank model. The pseudo-PageRank model outputs the nonnegative solution, which is normalized to obtain the stochastic solution if necessary. Pseudo-PageRank comes from a graph Laplacian matrix directly and is more suitable for practical applications \cite{LM06,Gl'15}.
% researchers replaced the columnwise-stochastic matrix in PageRank with a columnwise-substochastic matrix.
% a more general columnwise-substochastic matrix is a substitute for the columnwise stochastic matrix in PageRank. The resulting relaxation problem was called the pseudo-PageRank \cite{Gl'15}.
%This analogy allows us to compute the fixed point/dominant eigenvector of the PageRank eigensystem as the solution of the M-linear system, opening the door to a whole range of different numerical solvers for PageRank problems.
Moreover, the pseudo-PageRank model is potentially true of the multi-linear cases. This is the primary motivation of this paper.

% The connection between PR and pseudo-PR is a well known fact in the linear (matrix) case. I think this should be made clear. See for example [1-3, 4 section sign3.1]. In the matrix case, this analogy allows us to compute the fixed point/dominant eigenvector as the solution of a linear system, opening the door to a whole range of different numerical solvers. The same is potentially true in the multi-linear case, where the problem can be transformed into a minimization problem, as the authors show.

To explore higher-order extensions of classical Markov chain and PageRank for large-scale practical problems, on one hand, Li and Ng \cite{LN'14} used a rank-one array to approximate the limiting probability distribution array of a higher-order Markov chain. The vector for generating the rank-one array is called the limiting probability distribution vector of the higher-order Markov chain. On the other hand, Benson et al. \cite{BGL'17} studied the spacey random surfer, which remembers bits and pieces of history and is influenced by this information. Motivated by the spacey random surfer, Gleich et al. \cite{GLY'15} developed the multi-linear PageRank (MPR) model, of which the stochastic solution is called a multi-linear PageRank vector. Indeed, the multi-linear PageRank vector could be viewed as a limiting probability distribution vector of a special higher-order Markov chain.

According to the Brouwer fixed-point theorem, a stochastic solution of MPR always exists \cite{GLY'15}. The uniqueness of MPR solutions and related problems are more complicated than that of the classical PageRank. Many researchers gave various sufficient conditions \cite{CZ'13,FT'20,GLY'15,LLNV'17,LLVX'20,LN'14}. Dynamical systems of MPR were studied in \cite{BGL'17}. The perturbation analysis and error bounds of MPR solutions were presented in \cite{LCN'13,LLVX'20}. The ergodicity coefficient to a stochastic tensor was addressed in \cite{FT'20}. Interestingly, Benson \cite{B'19} studied three eigenvector centralities corresponding to tensors arising from hypergraphs. Tensor eigenvectors determine MPR vectors of these stochastic processes and capture some important higher-order network substructures \cite{BGL'15,BGL'17}.

For solving MPR and associated higher-order Markon chains, Li and Ng \cite{LN'14} proposed a fixed-point iteration, which converges linearly when a sufficient condition for uniqueness holds. Gleich et al. \cite{GLY'15} studied five numerical methods: a fixed-point iteration, a shifted fixed-point iteration, a nonlinear inner-outer iteration, an inverse iteration, and a Newton iteration. Meini and Poloni \cite{MP'18} solved a sequence of PageRank subproblems for Perron vectors to approximate the MPR solution. Furthermore, Huang and Wu \cite{HW'21} analyzed this Perron-based algorithm. Liu et al. \cite{LLV'19} designed and analyzed several relaxation algorithms for solving MPR. Cipolla et al. \cite{CRT'20} introduced extrapolation-based accelerations for a shifted fixed-point method and an inner-outer iteration. Hence, it is an active research area to compute MPR vectors arising from practical problems.

\subsection{Contributions}

In the process of extending the pseudo-PageRank method from graphs to uniform hypergraphs, the contributions are given as follows.

\begin{itemize}

  \item Using a Laplacian tensor of a uniform hypergraph directly without any dangling corrections, we propose the multi-linear pseudo-PageRank (MLPPR) model for a uniform hypergraph, which is a multi-linear system with nonnegative constraints. Every frontal slice of the coefficient tensor of MLPPR is an M-matrix. Comparing MLPPR with MPR, MLPPR admits a columuwise-substochastic tensor with dangling column fibers, while MPR rejects any dangling column fibers.

  \item The theoretical framework of the proposed model is established. The existence and uniqueness of the solution of MLPPR are discussed. The perturbation bound of the MLPPR solution is also presented.

  \item According to the structure of the tensor in MLPPR, we design a tensor splitting algorithm for finding a solution of MLPPR. Under suitable conditions, the propose algorithm converges linearly to a solution of MLPPR.

  \item As applications, we employ MLPPR to conduct hypergraph partitioning problems. For the problem of subspace clustering, we illustrate that MLPPR finds vertices of the target cluster more accurately when compared with some existing methods. For network analysis, directed 3-cycles in real-world networks are carefully selected for constructing edges of $3$-uniform hypergraphs. For such hypergraphs, numerical experiments illustrate that the new method is effective and efficient.
\end{itemize}

\subsection{Organization}

The outline of this paper is drawn as follows. Section 2 gives preliminary on tensors and graph-based PageRank. Hypergraphs and MLPPR are introduced in Section 3 and theoretical results on MLPPR are presented here. In Section 4, a tensor splitting algorithm for solving MLPPR is customized and analyzed. Applications in hypergraph partitioning problems are reported in Sections 5. Finally, we conclude this paper in Section 6.

\section{Preliminary}

Before starting, we introduce our notations and definitions.
%Let $\Re$ be the field of real numbers. We denote the vector space of all real $n$-dimensional vectors by $\Re^n$. Its nonnegative orthant is written as $\Re^n_+:=\{\xv\in\Re^n:x_i\ge0,i=1,2,\dots,n\}$. The notation $\Re^{m\times n}$ stands for the space of all $m$-by-$n$ real matrices.
We denote scalars, vectors, matrices, and tensors by lower case ($a$), bold lower case ($\av$), capital case ($A$), and calligraphic script ($\At$), respectively.
Let $\Re^{[k,n]}$ be the space of $k$th order $n$ dimensional real tensors and let $\Re^{[k,n]}_+$ be the set of nonnegative tensors in $\Re^{[k,n]}$. A tensor $\Pt\in\Re^{[k,n]}$ is called semi-symmetric if all of its elements $p_{i_1i_2\dots i_k}$ are invariant for any permutations of the second to the last indices $i_2,\dots,i_k$.

For $s\in\{1,\dots,k\}$, the $s$-mode product of $\Pt\in\Re^{[k,n]}$ and $\xv\in\Re^n$ produces a $(k-1)$th order $n$ dimensional tensor $\Pt\bar{\times}_s\xv$ of which elements are defined by
\begin{equation*}
  [\Pt\bar{\times}_s\xv]_{i_1\dots i_{s-1}i_{s+1}\dots i_k} := \sum_{i_s=1}^n p_{i_1\dots i_s\dots i_k}x_{i_s}.
\end{equation*}
For convenience, we define $\Pt\xv^{k-s}:=\Pt\bar{\times}_{s+1}\xv\cdots\bar{\times}_k\xv$. Particularly, $\Pt\xv^{k-1}$ and $\Pt\xv^{k-2}$ are a vector and a matrix, respectively.
For a $k$th order tensor $\Xt$ and an $\ell$th order tensor $\Yt$, the tensor outer product $\Xt\circ\Yt$ is a $(k+\ell)$th order tensor, of which elements are
\begin{equation*}
  [\Xt\circ\Yt]_{i_1\dots i_k i_{k+1}\dots i_{k+\ell}} = [\Xt]_{i_1\dots i_k}[\Yt]_{i_{k+1}\dots i_{k+\ell}}.
\end{equation*}
% for all $i_1,\dots,i_{k+\ell}$.
Particularly, when $k=\ell=1$, the outer product of two vectors $\xv$ and $\yv$ is a rank-one matrix $\xv\circ\yv=\xv\yv^T$.

Let $\one_n:=(1,1,\dots,1)^T\in\Re^n$ be an all-one vector. For convenience, we denote $\one=\one_n$. If $\vv\in\Re^n_+$ satisfies $\one^T\vv=1$, $\vv$ is called a stochastic vector. The $1$-mode unfolding of a tensor $\Pt\in\Re^{[k,n]}$ is an $n$-by-$n^{k-1}$ matrix $\Rm(\Pt)$ of which elements are
\begin{equation*}
  [\Rm(\Pt)]_{i\ell} := p_{ij_2j_3\dots j_k}, \quad\text{ where }\quad
  \ell=j_2+(j_3-1)n+\cdots+(j_k-1)n^{k-2}.
\end{equation*}
Each column of $\Rm(\Pt)$ is called a mode-1 (column) fiber of the tensor $\Pt$. If all of column fibers of $\Pt$ are stochastic, i.e., $\Pt\in\Re^{[k,n]}_+$ and $\one^T\Rm(\Pt)=\one^T_{n^{k-1}}$, $\Pt$ is called a columnwise-stochastic tensor. $\Pt$ is said to be a columnwise-substochastic tensor if $\Pt\in\Re^{[k,n]}_+$ and $\one^T\Rm(\Pt)\le\one^T_{n^{k-1}}$. Obviously, a columnwise-stochastic tensor is columnwise-substochastic and not conversely.
%Obviously, the columnwise-substochastic tensor admits zero fibers but the columnwise-stochastic tensor excludes zero fibers.
%The multistochastic tensor, of which all mode-$j$ fibers for $j=1,\dots,k$ are stochastic \cite{CLN14}, is a special case of both columnwise-stochastic tensors and columnwise-substochastic tensors.
In the remainder of this paper, a fiber refers to a mode-1 (column) fiber. For given indices $i_3,\dots,i_k$, an $n\times n$ matrix $(p_{iji_3\dots i_k})_{i,j=1,\dots,n}$ is called a frontal slice of a tensor $\Pt\in\Re^{[k,n]}$.

%Two useful products of a tensor $\Pt\in\Re^{[k,n]}$ and a vector $\xv\in\Re^n$ generate a vector $\Pt\xv^{k-1}$ and a matrix $\Pt\xv^{k-2}$ of which elements are defined respectively as
%\begin{eqnarray*}
%  {[\Pt\xv^{k-1}]}_i &=& \sum_{j_2=1}^n\cdots\sum_{j_k=1}^n P_{ij_2\dots j_k}x_{j_2}\cdots x_{j_k}, \\
%  {[\Pt\xv^{k-2}]}_{i_1i_2} &=& \sum_{j_3=1}^n\cdots\sum_{j_k=1}^n P_{i_1i_2j_3\dots j_k}x_{j_3}\cdots x_{j_k}.
%\end{eqnarray*}

% \vspace{100pt}

Next, we briefly review PageRank for directed graphs. A directed graph $G=(V,\vec{E},\wv)$ contains a vertex set $V:=\{1,2,\dots,n\}$, an arc set $\vec{E}:=\{\vec{e}_1,\vec{e}_2,\dots,\vec{e}_m\}$, and a weight vector $\wv\in\Re^m_+$. A component $w_{\ell}$ of $\wv$ stands for the weight of an arc $\vec{e}_{\ell}$ for $\ell=1,\dots,m$. An arc $\vec{e}=(j,i)\in\vec{E}$ means an directed edge $(j,i)$ from vertex $j$ to vertex $i$, where $j$ has an out-neighbor $i$ and $i$ has an in-neighbor $j$. The out-degree $d^{\mathrm{out}}_j$ of $j\in V$ is defined by the sum of weights of arcs that depart from $j$:
\begin{equation*}
  d^{\mathrm{out}}_j := \sum_{i: (j,i)=\vec{e}_{\ell}\in\vec{E}}w_{\ell} \ge0.
\end{equation*}
If $d^{\mathrm{out}}_j=0$, vertex $j$ is called a dangling vertex. Let the out-degree vector be  $\dv:=(d^{\mathrm{out}}_j)\in\Re^n_+$. The adjacency matrix $A=[a_{ij}]\in\Re^{n\times n}$ of $G$ is defined with elements
\begin{equation*}
  a_{ij} := \left\{\begin{aligned}
    & w_{\ell}, && \text{ if }(j,i)=\vec{e}_{\ell}\in\vec{E}, \\
    & 0,        && \text{ otherwise.}
  \end{aligned}\right.
\end{equation*}
Clearly, it holds that $\dv^T=\one^TA$. Vertex $j\in V$ is dangling if and only if the $j$th column of the adjacency matrix $A$ is a dangling (zero) vector.

Let us consider the random walk (stochastic process) $\{X(t):t=0,1,2,\dots\}$ on the vertex set $V=\{1,2,\dots,n\}$ of a directed graph $G$. The random walk takes a step according to the transition probability (columnwise-stochastic) matrix $P=(p_{ij})\in\Re^{n\times n}_+$ with probability $\alpha\in[0,1)$ and jumps randomly to a fixed distribution $\vv=(v_i)\in\Re^n_+$ with probability $1-\alpha$, i.e.,
\begin{equation}\label{prob}
  \prob\{X(t+1) = i \mid X(t) = j\} = \alpha p_{ij} + (1-\alpha)v_i,
\end{equation}
where $p_{ij}$ denotes the probability of a walk from vertex $j$ to vertex $i$. Here, $p_{ij}>0$ if and only if $(j,i)$ is an arc of $G$. % Also, $\sum_ip_{ij}=1$ for all $j$.
If the out-degree vector $\dv$ is positive, the normalized adjacency matrix $AD^+$ of $G$ is columnwise-stochastic, where $D:=\mathrm{diag}(\dv)$ and $(\cdot)^+$ stands for the Moore-Penrose pseudo-inverse. However, when the graph contains dangling vertices, the normalized adjacency matrix $AD^+$ is columnwise-substochastic but not columnwise-stochastic. Researchers \cite{Gl'15,LM06} incorporated a rank-one correction (named a dangling correction) into the columnwise-substochastic matrix to produce a columnwise-stochastic matrix:
\begin{equation*}
  P = AD^+ + \uv\cv^T,
\end{equation*}
where $\uv\in\Re^n$ is a stochastic vector, $\cv^T:=\one^T-\one^TAD^+$ is an indicator vector of dangling columns of $A$ and thus $\cv\in\{0,1\}^n$ is an indicator vector of dangling vertices in $V$. Here,  $c_j=1$ if the $j$th column of $A$ is a zero vector and $j\in V$ is a dangling vertex, $c_j=0$ otherwise.

Let a stochastic vector $\xv$ be a limiting probability distribution of the random walk $\{X(t)\}$. The random walk \eqref{prob} could be represented as a classical PageRank:
\begin{equation}\label{pr-2}
  \xv = \alpha P\xv + (1-\alpha)\vv,\quad \one^T\xv=1, \quad \xv\in\Re^n_+.
\end{equation}
It is well-known that the solution $\xv$ of PageRank \eqref{pr-2} could be interpreted as the Perron vector (the eigenvector corresponding to the dominant eigenvalue) of a nonnegative matrix, $(\alpha P + (1-\alpha)\vv\one^T)\xv=\xv$,
%\begin{equation*}
%  \xv = (\alpha P + (1-\alpha)\vv\one^T)\xv
%\end{equation*}
or the solution of a linear system with its coefficient matrix being an M-matrix, $(I-\alpha P)\xv = (1-\alpha)\vv$, where $I\in\Re^{n\times n}$ is the identity matrix.
%\Chg{Here a linear system $C\xv=\bv$ is called an M-linear system if its coefficient matrix $C$ is an M-matrix.}
%\begin{equation}\label{pr-linsys}
%  (I-\alpha P)\xv = (1-\alpha)\vv.
%\end{equation}
The linear system $(I-\alpha P)\xv = (1-\alpha)\vv$ with $\uv=\vv$ and $P = AD^+ + \uv\cv^T$ for the graph PageRank could be simplified by finding the solution $\yv\in\Re^n_+$ of another structural linear system
\begin{equation}\label{ppr}
  (I-\alpha AD^+)\yv=\vv
\end{equation}
and then normalizing solution $\yv$ to obtain a stochastic vector $\xv=\yv/(\one^T\yv)$ \cite{Gl'15}. In fact, the system \eqref{ppr} is called a pseudo-PageRank problem. The coefficient matrix $I-\alpha AD^+$ is both an M-matrix and a Laplacian matrix of the graph \cite{YCO'21}, enjoying good properties from linear algebra and graph theory. It is interesting to note that the pseudo-PageRank model doesn't need to do any dangling corrections. Indeed, dangling cases are inevitable for hypergraphs.

\section{Hypergraph and Multi-linear pseudo-PageRank}

Spectral hypergraph theory made great progress in the past decade \cite{BP'13,CD'12,QL17book,GCH22}. Let $G=(V,E,\wv)$ be a weighted undirected $k$-uniform hypergraph, where $V:=\{1,2,\dots,n\}$ is the vertex set, $E:=\{e_j\subseteq V:~\lvert e_j\rvert=k\text{ for }j=1,2,\dots,m\}$ is the edge set, and $\wv = (w_j)\in\Re^m_+$ is a positive vector whose component $w_j$ denotes the weight of an edge $e_j\in E$. Here, $\lvert\cdot\rvert$ means the cardinality of a set. Cooper and Dutle \cite{CD'12} defined adjacency tensors of uniform hypergraphs.

\begin{Def}[\cite{CD'12}]\label{AdjT-undir}
  Let $G=(V,E,\wv)$ be a weighted undirected $k$-uniform hypergraph with $n$ vertices. An adjacency tensor $\At=(a_{i_1\dots i_k})\in \Re^{[k,n]}_+$ of $G$ is a symmetric and nonnegative tensor of which elements are
  \begin{equation}\label{AdjTen-undir}
    a_{i_1\dots i_k}:=\left\{\begin{aligned}
      & \frac{w_j}{(k-1)!}, && \text{ if }\{i_1,\dots,i_k\}=e_j\in E, \\
      & 0,                  && \text{ otherwise. }
    \end{aligned}\right.
  \end{equation}
\end{Def}

Our MLPPR model could be applied in a weighted directed $k$-uniform hypergraph $G=(V,\vec{E},\wv)$ with a vertex set $V=\{1,2,\dots,n\}$, an arc set $\vec{E}:=\{\vec{e}_1,\vec{e}_2,\dots,\vec{e}_m\}$, and weights of arcs $\wv\in\Re^m_+$. An arc $\vec{e}_j\in\vec{E}$ is an ordered subset $\vec{e}_j:=(i_{j,1},i_{j,2},\dots,$ $i_{j,k})$ of $V$ for $j=1,2,\dots,m$. Here, vertices $i_{j,1},\dots,i_{j,k-1}$ are called tails of the arc $\vec{e}_j$ and the order among them are irrelevant. The last vertex $i_{j,k}$ is coined the head of the arc $\vec{e}_j$. Each arc $\vec{e}_j$ has a positive weight $w_j$ for $j=1,2,\dots,m$. Similarly, we define an adjacency tensor for the directed $k$-uniform hypergraph $G$ as follows.

\begin{Def}\label{AdjT-dir}
  Let $G=(V,\vec{E},\wv)$ be a weighted directed $k$-uniform hypergraph with $n$ vertices. An adjacency tensor $\At=(a_{i_1\dots i_k})\in \Re^{[k,n]}_+$ of $G$ is a nonnegative tensor of which elements are
  \begin{equation}\label{AdjTen-dir}
    a_{i_1\dots i_k}:=\left\{\begin{aligned}
      & \frac{w_j}{(k-1)!}, && \text{ if }(\sigma(i_2,\dots,i_k),i_1)=\vec{e}_j\in \vec{E}, \\
      & 0,                  && \text{ otherwise, }
    \end{aligned}\right.
  \end{equation}
  where $\sigma(i_2,\dots,i_k)$ denotes any permutation of indices $i_2,\dots,i_k$.
\end{Def}

Clearly, an adjacency tensor $\At$ of $G$ is semi-symmetric. We note that an edge $\{i_1,\dots,i_k\}$ of an undirected hypergraph could be viewed as $k$ varieties of arcs: $(\sigma(i_1,\dots,i_{k-1}),i_k), (\sigma(i_k,i_1,\dots,i_{k-2}),i_{k-1}), \dots, (\sigma(i_2,\dots,i_k),i_1)$ of a directed hypergraph. In this sense, Definitions \ref{AdjT-dir} and \ref{AdjT-undir} are consistent. % from the viewpoint of a vertex, the connection of the directed hypergraph is multiple-input and single-output. % This MISO structure is closely related to the following multi-linear pseudo-PageRank.

Next, we consider dangling fibers of an adjacency tensor arising from a uniform hypergraph. There are two kinds of dangling cases. (i) The structural dangling fibers are due to the definition of adjacency tensors of uniform hypergraphs. For example, since the 3-uniform hypergraph does not have edges/arcs with repeating vertices, elements $a_{ijj}$ of the associated adjacency tensor $\At$ are zeros for all $i$ and $j$. Hence, there are structural dangling fibers in adjacency tensors of any $k$-uniform hypergraphs with $k\ge3$. We note that a graph has no structural dangling fibers. (ii) The sparse dangling fibers are generated by the sparsity of a certain uniform hypergraph. For a set of tail vertices $\{i_1,\dots,i_{k-1}\}\subset V$, if there is no head vertex $i_k\in V$ such that  $(\sigma(i_1,\dots,i_{k-1}),i_k)\in \vec{E}$, there exist sparse dangling fibers. Particularly, there are amount of dangling fibers in adjacency tensors of large-scale sparse hypergraphs. A complete hypergraph does not have any sparse dangling fibers but has structural dangling fibers. % By the way, for large-scale sparse hypergraphs arising from real-world network datasets, there are amount of dangling fibers in the corresponding adjacency tensors.

\begin{Exm}
  Let $G=(V,E,\wv)$ be a $3$-uniform undirected hypergraph with $V=\{1,2,3,4\}$, $E=\{\{1,2,3\},\{2,3,4\}\}$, and $\wv=\one_4$. By Definition \ref{AdjT-undir}, the unfolding $\Rm(\At)\in\Re^{4\times 16}$ of the adjacency tensor $\At\in\Re^{[3,4]}$ of $G$ is written as follows:
  \begin{equation*}
  \bordermatrix{
    & \text{\scriptsize 11} & \text{\scriptsize 21} & \text{\scriptsize 31} & \text{\scriptsize 41} & \text{\scriptsize 12} & \text{\scriptsize 22} & \text{\scriptsize 32} & \text{\scriptsize 42} & \text{\scriptsize 13} & \text{\scriptsize 23} & \text{\scriptsize 33} & \text{\scriptsize 43} & \text{\scriptsize 14} & \text{\scriptsize 24} & \text{\scriptsize 34} & \text{\scriptsize 44} \cr
  \text{\scriptsize 1} &    &    &    &    &    &    &  \frac{1}{2} &    &    &  \frac{1}{2} &    &    &    &    &    &    \cr
  \text{\scriptsize 2} &    &    &  \frac{1}{2} &    &    &    &    &    &  \frac{1}{2} &    &    &  \frac{1}{2} &    &    &  \frac{1}{2} &    \cr
  \text{\scriptsize 3} &    & \frac{1}{2} &    &    &  \frac{1}{2} &    &    &  \frac{1}{2} &    &    &    &    &    &  \frac{1}{2} &    &    \cr
  \text{\scriptsize 4} &    &    &    &    &    &    &  \frac{1}{2} &    &    &  \frac{1}{2} &    &    &    &    &    &    \cr  }.
  \end{equation*}
  Obviously, columns $11$, $22$, $33$, and $44$ are structural dangling fibers. Columns $41$ and $14$ are sparse dangling fibers.
\end{Exm}

To get rid of these dangling fibers that lead to a gap between Laplacian tensors of uniform hypergraphs and columnwise-stochastic tensors for higher-order extensions of PageRank, we are going to develop a novel MLPPR model. Following the spirit of the pseudo-PageRank model for graphs, we normalize all nonzero fibers of the adjacency tensor $\At$ of a uniform hypergraph $G$ to produce a columnwise-substochastic tensor $\bar{\Pt}=(\bar{p}_{i_1\dots i_k})\in\Re^{[k,n]}_+$, of which elements are defined by
\begin{equation}\label{StoT}
  \bar{p}_{i_1\dots i_k} = \left\{\begin{aligned}
    & \frac{a_{i_1\dots i_k}}{\sum\limits_{\ell=1}^n a_{\ell i_2\dots i_k}},
      && \text{ if }\sum_{\ell=1}^n a_{\ell i_2\dots i_k}>0, \\
    & 0, && \text{ otherwise.}
    \end{aligned}\right.
\end{equation}
That is to say, we have $\Rm(\bar{\Pt})=\Rm(\At)D^+\in\Re^{n\times n^{k-1}}_+$, where $D:=\mathrm{diag}(\one^T\Rm(\At))$.
Both an adjacency tensor $\At$ and its normalization $\bar{\Pt}$ are nonnegative tensors, which have interesting theoretical results \cite{LN15}.

As an analogue of the Laplacian matrix of a graph \cite{YCO'21}, we call
\begin{equation*}
  I\circ\one^{\circ(k-2)}-\alpha\bar{\Pt}\in\Re^{[k,n]}
\end{equation*}
a Laplacian tensor of a $k$-uniform hypergraph $G$ with $n$ vertices, where $\alpha\in[0,1)$ is a probability and $\one^{\circ(k-2)}$ is the outer produce of $k-2$ copies of $\one$. Obviously, each frontal slice of this Laplacian tensor is an M-matrix. Now, we are ready to present the multi-linear pseudo-PageRank model.

\begin{Def}[MLPPR]\label{mpPR-def}
  Let $\bar{\Pt}$ be a columnwise-substochastic tensor, i.e., $\bar{\Pt}\in\Re^{[k,n]}_+$ and $\one^T\Rm(\bar{\Pt})\le\one^T_{n^{k-1}}$. Let $\vv\in\Re^n$ be a stochastic vector and let $\alpha\in[0,1)$ be a probability. Then MLPPR is to find a nonnegative solution $\yv\in\Re^n_+$ of the following multi-linear system
  \begin{equation}\label{mlPPR}
    \left(I\circ\one^{\circ(k-2)}-\alpha\bar{\Pt}\right)\yv^{k-1}=\vv.
  \end{equation}
  The normalized $\yv$, i.e., $\yv/(\yv^T\one)$, is called an MLPPR vector.
\end{Def}

%The tensor $I\circ\one\circ\cdots\circ\one-\alpha\bar{\Pt}$ in \eqref{mlPPR} can be interpreted as a generalization of a type of graph Laplacian for hypergraphs \cite{YCO'21}.

Next, we show some properties of the multi-linear system of MLPPR.

\begin{Lem}\label{Lem:3-4}
  % Let $\bar{\Pt}$, $\vv$, and $\alpha$ be defined as in Definition \ref{mpPR-def}.
  If the MLPPR system \eqref{mlPPR} has a nonnegative solution $\yv_*$, then $\yv_*$ belongs to a bounded, closed and convex set:
  \begin{equation}\label{y-bound'}
    \Delta:=\{\yv\in\Re^n_+ : \yv^T\one \le \varrho\},
  \end{equation}
  where $\varrho:= (1-\alpha)^{-\frac{1}{k-1}}$.
\end{Lem}
\begin{proof}
  Since $\bar{\Pt}$ is columnwise-substochastic, by taking inner products between $\one$ and both sides of \eqref{mlPPR}, we have
  \begin{equation}\label{y-bound'-1}
    1=\one^T\vv = (\one^T\yv_*)^{k-1}-\alpha \one^T(\bar{\Pt}\yv_*^{k-1})
      \ge (1-\alpha)(\one^T\yv_*)^{k-1},
  \end{equation}
  which means that $\yv_*$ satisfies $\one^T\yv_* \le (1-\alpha)^{-\frac{1}{k-1}}=\varrho$.
  %\begin{equation*}
  %  \one^T\yv_* \le (1-\alpha)^{-\frac{1}{k-1}}=\varrho.
  %\end{equation*}
  Hence, $\yv_*\in\Delta$. It is straightforward to see that $\Delta$ is a bounded, closed and convex set.
\end{proof}

\begin{Lem}\label{Lem:3-5}
  % Let $\bar{\Pt}$, $\vv$, and $\alpha$ be defined as in Definition \ref{mpPR-def}.
  Let $\yv_*\in\Delta$. Then $\yv_*$ solves the MLPPR system \eqref{mlPPR} if and only if $\yv_*$ is a fixed point of the continuous map
  \begin{equation}\label{map-PHI}
    \Phi(\yv) := \left(1+\alpha\one^T(\bar{\Pt}\yv^{k-1})\right)^{-\frac{k-2}{k-1}}(\vv+\alpha\bar{\Pt}\yv^{k-1}).
  \end{equation}
\end{Lem}
\begin{proof}
  On one hand, assume that $\yv_*$ solves MLPPR. From \eqref{y-bound'-1}, we say
  \begin{equation*}
    \one^T\yv_* = \left(1+\alpha\one^T(\bar{\Pt}\yv_*^{k-1})\right)^{\frac{1}{k-1}}.
  \end{equation*}
  Then the MLPPR system \eqref{mlPPR} could be represented as
  \begin{equation*}
    \yv_* = (\one^T\yv_*)^{-(k-2)}(\vv+\alpha\bar{\Pt}\yv_*^{k-1}) = \left(1+\alpha\one^T(\bar{\Pt}\yv_*^{k-1})\right)^{-\frac{k-2}{k-1}}(\vv+\alpha\bar{\Pt}\yv_*^{k-1}).
  \end{equation*}
  That is to say, $\yv_*$ is a fixed point of $\Phi$.

  On the other hand, From $\yv_*=\Phi(\yv_*)$ we can derive
  \begin{equation*}
    \one^T\yv_* = \one^T\Phi(\yv_*) = \left(1+\alpha\one^T(\bar{\Pt}\yv_*^{k-1})\right)^{\frac{1}{k-1}}.
  \end{equation*}
  Then the equality \eqref{map-PHI} reduces to
  \begin{equation*}
    \yv_* = (\one^T\yv_*)^{-(k-2)}(\vv+\alpha\bar{\Pt}\yv_*^{k-1}),
  \end{equation*}
  which is equivalent to MLPPR \eqref{mlPPR}.
\end{proof}

%The multi-linear system \eqref{mlPPR} exists a nonnegative solution as shown in the following theorem.
The following theorem reveals the existence of nonnegative solutions of MLPPR.

\begin{Thm}\label{Thm:3-6}
  Let $\bar{\Pt}$, $\vv$, and $\alpha$ be defined as in Definition \ref{mpPR-def}.
  Then the MLPPR system \eqref{mlPPR} has at least one nonzero solution $\yv_*\in\Re^n_+$.
\end{Thm}
\begin{proof}
  At the beginning, we prove that $\Phi$ maps $\Delta$ into itself, where $\Delta$ and $\Phi$ are defined as in Lemmas \ref{Lem:3-4} and \ref{Lem:3-5}, respectively. Let $\yv\in\Delta$. Obviously, $\Phi(\yv)\in\Re^n_+$. By taking an inner product between $\one$ and $\Phi(\yv)$, we have
  \begin{eqnarray*}
    \one^T\Phi(\yv) %&=& \left(1+\alpha\one^T(\bar{\Pt}\yv^{k-1})\right)^{-\frac{k-2}{k-1}+1} \\
      &=& \left(1+\alpha\one^T(\bar{\Pt}\yv^{k-1})\right)^{\frac{1}{k-1}} \\
      &\le& \left(1+\alpha(\one^T\yv)^{k-1}\right)^{\frac{1}{k-1}} \\
      &\le& \left(1+\alpha(1-\alpha)^{-1}\right)^{\frac{1}{k-1}} \\
      &=& \varrho,
  \end{eqnarray*}
  %\begin{equation*}
  %  \one^T\Phi(\yv) = \left(1+\alpha\one^T(\bar{\Pt}\yv^{k-1})\right)^{\frac{1}{k-1}}
  %    \le \left(1+\alpha(\one^T\yv)^{k-1}\right)^{\frac{1}{k-1}}
  %    \le \left(1+\alpha(1-\alpha)^{-1}\right)^{\frac{1}{k-1}}
  %    =\varrho,
  %\end{equation*}
  where the first inequality holds because $\bar{\Pt}$ is columnwise-substochastic.
  Hence, the continuous function $\Phi$ maps the compact and convex set $\Delta$ into itself.
  %we say $\Phi:\Delta\to\Delta$.

  According to the Brouwer fixed point theorem, there exists $\yv_*\in\Delta$ such that $\Phi(\yv_*)=\yv_*$.
  By Lemma \ref{Lem:3-5}, $\yv_*$ is also a nonnegative solution of MLPPR.
\end{proof}

\begin{Cor}\label{Cor:3-7}
  % Let $\bar{\Pt}$, $\vv$, and $\alpha$ be defined as in Definition \ref{mpPR-def}.
  If $\vv$ is a positive (stochastic) vector, then a nonnegative solution $\yv_*=\Phi(\yv_*)$ of MLPPR is also positive.
\end{Cor}

On the uniqueness of MLPPR solution, we have theoretical results under some conditions. To move on, the following two lemmas are helpful.

\begin{Lem}\label{Lem:3-8}
  Let $\otimes$ denote the Kronecker product. Define
  \begin{equation*}
    \xv^{\otimes(k-1)} := \underbrace{\xv\otimes\xv\otimes\cdots\otimes\xv}_{k-1\text{ copies of }\xv}.
  \end{equation*}
  For any $\xv,\yv\in\Delta$, we have
  \begin{equation*}
    \|\xv^{\otimes(k-1)}-\yv^{\otimes(k-1)}\|_1 \le (k-1)(1-\alpha)^{-\frac{k-2}{k-1}}\|\xv-\yv\|_1.
  \end{equation*}
\end{Lem}
\begin{proof}
  By direct calculations with \eqref{y-bound'}, we have
  \begin{eqnarray*}
    \lefteqn{ \|\xv^{\otimes(k-1)}-\yv^{\otimes(k-1)}\|_1 }\\
      &=& \|(\xv-\yv)\otimes\xv^{\otimes(k-2)}+\yv\otimes(\xv-\yv)\otimes\xv^{\otimes(k-3)}+\cdots+\yv^{\otimes(k-2)}\otimes(\xv-\yv)\|_1 \\
      &\le& \|\xv-\yv\|_1(\one^T\xv)^{k-2}+(\one^T\yv)\|\xv-\yv\|_1(\one^T\xv)^{k-3}+\cdots+(\one^T\yv)^{k-2}\|\xv-\yv\|_1  \\
      &\le& (k-1) (1-\alpha)^{-\frac{k-2}{k-1}} \|\xv-\yv\|_1,
  \end{eqnarray*}
  where the last inequality is valid because $\xv,\yv\in\Delta$.
\end{proof}

\begin{Lem}\label{Lem:3-9}
  Define $\varphi(\zv):=(\one^T\zv)^{-\frac{k-2}{k-1}}\zv$ for $\zv\in\Xi:=\{\zv\in\Re^n_+:\one^T\zv\ge1\}$. For any $\av,\bv\in\Xi$, we have
  \begin{equation*}
    \|\varphi(\av)-\varphi(\bv)\|_1 \le \frac{2k-3}{k-1}\|\av-\bv\|_1.
  \end{equation*}
\end{Lem}
\begin{proof}
  Let $\psi(t):=t^{-\frac{k-2}{k-1}}$ for $t\in[1,\infty)$. Obviously, $\psi(t)$ is a convex and monotonically decreasing function on $[1,\infty)$. When $s \ge t \ge 1$, we have
  \begin{equation}\label{3-8-a}
    \phi(s)-\phi(t)\le \phi'(s)(s-t)
  \end{equation}
  and
  \begin{equation}\label{3-8-b}
    s^{-\frac{k-2}{k-1}} \le t^{-\frac{k-2}{k-1}} \le 1.
  \end{equation}

  For any $\av,\bv\in\Xi$, without loss of generality, we suppose $\one^T\av\ge\one^T\bv\ge1$. Then, we have
  \begin{eqnarray*}
    \|\varphi(\av)-\varphi(\bv)\|_1
      &=& \left\|(\one^T\av)^{-\frac{k-2}{k-1}}\av-(\one^T\bv)^{-\frac{k-2}{k-1}}\bv\right\|_1 \\
      % &=& \left\|(\one^T\av)^{-\frac{k-2}{k-1}}(\av-\bv) +\left[(\one^T\av)^{-\frac{k-2}{k-1}}-(\one^T\bv)^{-\frac{k-2}{k-1}}\right]\bv\right\|_1 \\
      &\le& (\one^T\av)^{-\frac{k-2}{k-1}}\|\av-\bv\|_1 + \left\lvert(\one^T\av)^{-\frac{k-2}{k-1}}-(\one^T\bv)^{-\frac{k-2}{k-1}}\right\rvert\|\bv\|_1 \\
      &=& (\one^T\av)^{-\frac{k-2}{k-1}}\|\av-\bv\|_1 + \left[(\one^T\bv)^{-\frac{k-2}{k-1}}-(\one^T\av)^{-\frac{k-2}{k-1}}\right](\one^T\bv) \\
      &\le& (\one^T\av)^{-\frac{k-2}{k-1}}\|\av-\bv\|_1 -\frac{k-2}{k-1}(\one^T\bv)^{-\frac{k-2}{k-1}-1}(\one^T\bv-\one^T\av)(\one^T\bv) \\
      %&\le& \|\av-\bv\|_1 + \frac{k-2}{k-1}(\theta\one^T\av+(1-\theta)\one^T\bv)^{-\frac{k-2}{k-1}-1}(\one^T\av-\one^T\bv)(\one^T\bv) \\
      % &=& (\one^T\av)^{-\frac{k-2}{k-1}}\|\av-\bv\|_1 + \frac{k-2}{k-1}(\one^T\bv)^{-\frac{k-2}{k-1}}\one^T(\av-\bv) \\
      &\le& \|\av-\bv\|_1 + \frac{k-2}{k-1}\one^T(\av-\bv) \\
      %&\le& \|\av-\bv\|_1 + \frac{k-2}{k-1}\one^T(\av-\bv) \\
      % &\le& \|\av-\bv\|_1 + \frac{k-2}{k-1}\|\one\|_{\infty}\|\av-\bv\|_1 \\
      &\le& \frac{2k-3}{k-1}\|\av-\bv\|_1,
  \end{eqnarray*}
  where the second and the third inequalities hold owing to \eqref{3-8-a} and \eqref{3-8-b}, respectively.
\end{proof}

\begin{Thm}\label{Thm:3-A}
  Let $\bar{\Pt}$, $\vv$, and $\alpha$ be defined as in Definition \ref{mpPR-def}. Define
  \begin{equation}\label{varsigma}
    \varsigma := (2k-3)\alpha(1-\alpha)^{-\frac{k-2}{k-1}}.
  \end{equation}
  If $\varsigma<1$, then the map $\Phi:\Delta\to\Delta$ is a contraction map, i.e., $\|\Phi(\xv)-\Phi(\yv)\|_1\le\varsigma\|\xv-\yv\|_1$ for all $\xv,\yv\in\Delta$.
\end{Thm}
\begin{proof}
  For any $\xv,\yv\in\Delta$, we have $\vv+\alpha\bar{\Pt}\xv^{k-1},\vv+\alpha\bar{\Pt}\yv^{k-1}\in\Xi$ and
  \begin{eqnarray*}
    \|\Phi(\xv)-\Phi(\yv)\|_1
      &=& \|\varphi(\vv+\alpha\bar{\Pt}\xv^{k-1})-\varphi(\vv+\alpha\bar{\Pt}\yv^{k-1})\|_1 \\
      &\le& \frac{2k-3}{k-1}\|\vv+\alpha\bar{\Pt}\xv^{k-1}-\vv-\alpha\bar{\Pt}\yv^{k-1}\|_1 \\
      %&=& \alpha\frac{2k-3}{k-1}\|\bar{\Pt}\xv^{k-1}-\bar{\Pt}\yv^{k-1}\|_1 \\
      &=& \alpha\frac{2k-3}{k-1}\|\Rm(\bar{\Pt})(\xv^{\otimes(k-1)}-\yv^{\otimes(k-1)})\|_1 \\
      %&\le& \alpha\frac{2k-3}{k-1}\|\Rm(\bar{\Pt})\|_{1}\|\xv^{\otimes(k-1)}-\yv^{\otimes(k-1)}\|_1 \\
      &\le& \alpha\frac{2k-3}{k-1}\|\xv^{\otimes(k-1)}-\yv^{\otimes(k-1)}\|_1 \\
      &\le& \alpha\frac{2k-3}{k-1}(k-1)(1-\alpha)^{-\frac{k-2}{k-1}}\|\xv-\yv\|_1 \\
      &=& \varsigma\|\xv-\yv\|_1,
  \end{eqnarray*}
  where the first and the last inequalities are owing to Lemmas \ref{Lem:3-9} and \ref{Lem:3-8}, respectively. The second inequality holds because $\bar{\Pt}$ is columnwise-substochastic. Hence, $\Phi$ is a contraction map when $\varsigma<1$.
\end{proof}

For the special case of the pseudo-PageRank model for graphs, we know $k=2$. Then the condition $\varsigma<1$ reduces to a trivial condition $\alpha<1$. When $k=3$, the condition $\varsigma<1$ reduces to $\alpha<\frac{\sqrt{37}-1}{18}\approx 0.28$. Furthermore, the uniqueness of MLPPR solution follows this contraction theorem.

\begin{Cor}\label{Cor:3-B}
  Under conditions of Theorem \ref{Thm:3-A}, the MLPPR system \eqref{mlPPR} has a unique nonnegative solution $\yv_*\in\Delta$.
\end{Cor}

%\begin{Cor}\label{Cor:3-C}
%  Under conditions of Theorem \ref{Thm:3-A}. We consider an algorithm:
%  \begin{equation*}
%    \text{ Choose } \yv_0\in\Delta, \text{ compute }\yv_{c+1}=\Phi(\yv_c)\text{ for }c=0,1,2,\dots.
%  \end{equation*}
%  Then the sequence $\{\yv_c\}$ converges to $\yv_*$ with a linear convergence rate.
%\end{Cor}

Whereafter, we will consider how perturbations in the columnwise-substochastic tensor $\bar{\Pt}$ and the stochastic vector $\vv$ affect the MLPPR solution $\yv$. Numerical verifications of this theorem will be illustrated in Subsection \ref{toy-hgraph}.

\begin{Thm}\label{Thm:3-D}
  Let $\bar{\Pt}$ and $\bar{\Pt}_p=\bar{\Pt}+\delta\bar{\Pt}\in\Re^{[k,n]}$ be columnwise-substochastic tensors. Let $\vv$ and $\vv_p=\vv+\delta\vv\in\Re^n$ be stochastic vectors. The probability $\alpha$ is chosen such that $\varsigma<1$. Nonnegative vectors $\yv\in\Delta$ and $\yv_p=\yv+\delta\yv\in\Delta$ solve MLPPR systems
  \begin{equation*}
    \left(I\circ\one^{\circ(k-2)}-\alpha\bar{\Pt}\right)\yv^{k-1}=\vv \quad\text{ and }\quad
    \left(I\circ\one^{\circ(k-2)}-\alpha\bar{\Pt}_p\right)\yv_p^{k-1}=\vv_p,
  \end{equation*}
  respectively. Then, we have
  \begin{eqnarray*}
    \|\delta\yv\|_1 &\le& \frac{2k-3}{(k-1)(1-\varsigma)}
    \left(\frac{\alpha}{1-\alpha}\|\Rm(\delta\bar{\Pt})\|_1+\|\delta\vv\|_1\right).
  \end{eqnarray*}
\end{Thm}
\begin{proof}
  According to definition of $\varphi$ in Lemma \ref{Lem:3-9}, we have
  \begin{eqnarray*}
    \lefteqn{ \|\delta\yv\|_1 = \|\yv_p-\yv\|_1 }\\
      &=& \|\varphi(\vv_p+\alpha\bar{\Pt}_p\yv_p^{k-1})-\varphi(\vv+\alpha\bar{\Pt}\yv^{k-1})\|_1 \\
      &\le& \frac{2k-3}{k-1}\|\vv_p+\alpha\bar{\Pt}_p\yv_p^{k-1}-\vv-\alpha\bar{\Pt}\yv^{k-1}\|_1 \\
      &=& \frac{2k-3}{k-1}\|\alpha(\bar{\Pt}_p-\bar{\Pt})\yv_p^{k-1} + \alpha(\bar{\Pt}\yv_p^{k-1}-\bar{\Pt}\yv^{k-1})+\delta\vv\|_1\\
      &\le& \alpha\frac{2k-3}{k-1}\|\Rm(\delta\bar{\Pt})\|_1\|\yv_p^{\otimes (k-1)}\|_1
      + \alpha\frac{2k-3}{k-1}\|\Rm(\bar{\Pt})\|_1\|\yv_p^{\otimes(k-1)}-\yv^{\otimes(k-1)}\|_1 \\&&{}
      + \frac{2k-3}{k-1}\|\delta\vv\|_1\\
      &\le& \alpha\frac{2k-3}{k-1}\|\Rm(\delta\bar{\Pt})\|_1(\one^T\yv_p)^{k-1} %\\&&{}
      + \alpha\frac{2k-3}{k-1}\|\yv_p^{\otimes(k-1)}-\yv^{\otimes(k-1)}\|_1
      + \frac{2k-3}{k-1}\|\delta\vv\|_1  \\
      &\le& \frac{\alpha}{1-\alpha}\frac{2k-3}{k-1}\|\Rm(\delta\bar{\Pt})\|_1
      + \alpha(1-\alpha)^{-\frac{k-2}{k-1}}(2k-3)\|\delta\yv\|_1 %\\&&{}
      + \frac{2k-3}{k-1}\|\delta\vv\|_1,
  \end{eqnarray*}
  where the first inequality is valid by Lemma \ref{Lem:3-9} and the last inequality follows from $\yv_d\in\Delta$ and Lemma \ref{Lem:3-8}. Hence, we obtain
  \begin{equation*}
    % \frac{1-\varsigma}{2k-3} \|\delta\yv\|_1     =
    \left(\frac{1}{2k-3}-\alpha(1-\alpha)^{-\frac{k-2}{k-1}}\right)\|\delta\yv\|_1
    \le \frac{1}{k-1}\left(\frac{\alpha}{1-\alpha}\|\Rm(\delta\bar{\Pt})\|_1+ \|\delta\vv\|_1\right).
  \end{equation*}
  The proof is completed.
\end{proof}

\subsection{Relationships between MPR and MLPPR}

In the framework of MPR \cite{BGL'15,GLY'15}, a columnwise-substochastic tensor $\bar{\Pt}$ used in MLPPR was adjusted to a columnwise-stochastic tensor by performing the following dangling correction:
\begin{equation}\label{dePseudo}
  \Pt := \bar{\Pt}+\vv\circ\left(\one^{\circ(k-1)}-\bar{\Pt}\bar{\times}_1 \one\right) \in\Re^{[k,n]}_+.
\end{equation}
We note that the tensor $\one^{\circ(k-1)}-\bar{\Pt}\bar{\times}_1\one \in\Re^{[k-1,n]}$ is dense if $\bar{\Pt}\bar{\times}_1\one$ is sparse. Gleich et al. \cite{GLY'15} proposed the MPR model:
\begin{equation}\label{mlPR}
  \xv = \alpha\Pt\xv^{k-1}+(1-\alpha)\vv, \quad \xv\in\Re^n_+, \quad \one^T\xv=1,
\end{equation}
where $\Pt\in\Re^{[k,n]}$ is a columnwise-stochastic tensor and $\vv\in\Re^n$ is a stochastic vector.

We remark that MPR \eqref{mlPR} has a homogeneous form
\begin{equation}\label{hmlPR}
  \zv(\one^T\zv)^{k-2} = \alpha\Pt\zv^{k-1}+(1-\alpha)\vv(\one^T\zv)^{k-1}, \quad \zv\in\Re^n_+.
\end{equation}
If the homogeneous form \eqref{hmlPR} has a nonzero solution $\zv$, there are an infinite number of solutions $\gamma\zv$ for all $\gamma>0$. Furthermore, a nonzero vector $\zv$ solves the homogeneous form \eqref{hmlPR} if and only if $\zv/(\one^T\zv)$ solves MPR \eqref{mlPR}.

Next, we discuss relationships between MLPPR, MPR, and the homogeneous MPR when the dangling correction \eqref{dePseudo} holds.

\begin{Lem}\label{Lem:3-E}
  % Suppose that $\bar{\Pt}$, $\vv$, and $\alpha$ are introduced in Definition \ref{mpPR-def} and \eqref{dePseudo} is valid.
  Suppose the dangling correction \eqref{dePseudo} holds.
  If $\xv$ solves MPR \eqref{mlPR}, then
  \begin{equation*}
    \yv:=\left(1-\alpha\one^T(\bar{\Pt}\xv^{k-1})\right)^{-\frac{1}{k-1}}\xv
  \end{equation*}
  is a nonnegative solution of MLPPR \eqref{mlPPR}.
\end{Lem}
\begin{proof}
  From \eqref{mlPR} and \eqref{dePseudo}, we have
  \begin{eqnarray*}
    \xv &=& \alpha\Pt\xv^{k-1}+(1-\alpha)\vv \\
      &=& \alpha\bar{\Pt}\xv^{k-1} + \alpha\vv - \alpha\vv\one^T(\bar{\Pt}\xv^{k-1}) + (1-\alpha)\vv \\
      &=& \alpha\bar{\Pt}\xv^{k-1} + \left(1-\alpha\one^T(\bar{\Pt}\xv^{k-1})\right)\vv,
  \end{eqnarray*}
  which implies
  \begin{equation*}
    (I\circ\one^{\circ(k-2)}-\alpha\bar{\Pt})\xv^{k-1}\left(1-\alpha\one^T(\bar{\Pt}\xv^{k-1})\right)^{-1}=\vv.
  \end{equation*}
  The proof is completed.
\end{proof}

\begin{Lem}\label{Lem:3-F}
  % Suppose that $\bar{\Pt}$, $\vv$, and $\alpha$ are introduced in Definition \ref{mpPR-def} and \eqref{dePseudo} is valid.
  Suppose the dangling correction \eqref{dePseudo} holds.
  If a nonnegative vector $\yv$ solves MLPPR \eqref{mlPPR}, then $\yv$ is a solution of the homogeneous form of MPR \eqref{hmlPR}.
\end{Lem}
\begin{proof}
  By taking inner products between $\one$ and both sides of \eqref{mlPPR}, we get
  \begin{equation*}
    \alpha\one^T(\bar{\Pt}\yv^{k-1})=(\one^T\yv)^{k-1}-1.
  \end{equation*}
  Then, it follows from \eqref{mlPPR} and \eqref{dePseudo} that
  \begin{eqnarray*}
    \vv &=& \yv(\one^T\yv)^{k-2}-\alpha\bar{\Pt}\yv^{k-1} \\
      &=& \yv(\one^T\yv)^{k-2}-\alpha\Pt\yv^{k-1} + \alpha\vv(\one^T\yv)^{k-1}-\alpha\vv\one^T(\bar{\Pt}\yv^{k-1})  \\
      &=& \yv(\one^T\yv)^{k-2}-\alpha\Pt\yv^{k-1} + \alpha\vv(\one^T\yv)^{k-1}-\vv((\one^T\yv)^{k-1}-1) \\
      &=& \yv(\one^T\yv)^{k-2}-\alpha\Pt\yv^{k-1} - (1-\alpha)\vv(\one^T\yv)^{k-1}+\vv,
  \end{eqnarray*}
  which implies \eqref{hmlPR} straightforwardly.
\end{proof}

According to Lemmas \ref{Lem:3-E} and \ref{Lem:3-F}, we get the following unexpected but reasonable theorem.

\begin{Thm}\label{Thm:3-G}
  Let $\bar{\Pt}$, $\vv$, and $\alpha$ be defined as in Definition \ref{mpPR-def}. Suppose the dangling correction \eqref{dePseudo} holds. Then the solution of MPR \eqref{mlPR} is unique if and only if the nonnegative solution of MLPPR \eqref{mlPPR} is unique.
\end{Thm}
\begin{proof}
  Suppose that $\yv_1$ and $\yv_2$ are nonnegative solutions of MLPPR \eqref{mlPPR}. By Lemma \ref{Lem:3-F}, $\yv_1/(\one^T\yv_1)$ and $\yv_2/(\one^T\yv_2)$ are solutions of MPR \eqref{mlPR}. Because MPR has a unique solution, we have
  \begin{equation*}
    \yv_1/(\one^T\yv_1) = \yv_2/(\one^T\yv_2),
  \end{equation*}
  which means $\yv_1=\gamma\yv_2$ for a positive number $\gamma$.

  On the other hand, because $\yv_1$ and $\yv_2$ solve MLPPR \eqref{mlPPR}, we have
  \begin{eqnarray*}
    \vv &=& (I\circ\one^{\circ(k-2)}-\alpha\bar{\Pt})\yv_1^{k-1} \\
      &=& (I\circ\one^{\circ(k-2)}-\alpha\bar{\Pt})(\gamma\yv_2)^{k-1} \\
      &=& \gamma^{k-1}(I\circ\one^{\circ(k-2)}-\alpha\bar{\Pt})\yv_2^{k-1} \\
      &=& \gamma^{k-1}\vv.
  \end{eqnarray*}
  Hence, the positive number $\gamma$ is subject to $\gamma^{k-1}-1=0$. We say $\gamma=1$ which means $\yv_1=\yv_2$.

  The converse statement could be verified by using Lemma \ref{Lem:3-E} and a similar discussion.
\end{proof}

Based on Theorem \ref{Thm:3-G} and Theorem 4.3 in \cite{GLY'15}, we get a new sufficient condition for preserving the uniqueness of an MLPPR solution.

\begin{Cor}\label{Cor:3-H}
  MLPPR has a unique nonnegative solution when $\alpha<\frac{1}{k-1}$.
\end{Cor}

Clearly, Corollary \ref{Cor:3-H} is better than Corollary \ref{Cor:3-B}, since we can infer $\alpha<\frac{1}{k-1}$ from $\varsigma=(2k-3)\alpha(1-\alpha)^{-\frac{k-2}{k-1}}<1$ and not conversely when $k\ge3$. In the case of $k=2$, $\alpha<\frac{1}{k-1}$ is equivalent to $\varsigma<1$.
Additionally, other sufficient conditions for the uniqueness of the MPR solution elaborated in \cite{LN'14,LLNV'17,LLVX'20} are also applied to MLPPR.

Comparing with MLPPR, MPR adjusts all fibers of the columnwise-substochastic tensor to stochastic vectors in advance and then MPR returns a stochastic solution. When applied into uniform hypergraphs, the adjustment in MPR lacks significance of graph theory. In addition, it is a time-consuming process to produce the adjustment for large scale hypergraphs. Our MLPPR utilizes the columnwise-substochastic tensor directly to produce a nonnegative solution, which could easily be normalized to a stochastic vector if necessary. By reversing the order of adjustment/normalization and root-finding of constrained multi-linear systems, we proposed a novel MLPPR model.
Moreover, MLPPR gets rid of dangling cases by directly dealing with columnwise-substochastic tensors, which are usually sparser than columnwise-stochastic tensors used in MPR. This sparsity reveals a valuable potential of MLPPR in computation.

% On the other hand, the domain of MLPPR $\Re^n_+$ is simpler than that of MPR $\{\xv\in\Re^n_+:\one^T\xv=1\}$. This enables us to customize an easy and effective numerical algorithm for solving MLPPR.
% In fact, the positive orthant $\Re^n_+$ often arose in tensor problems \cite{LGL16,HST19}.

\section{Tensor splitting algorithm}

Inspired by the structure of the coefficient tensor of MLPPR, of which all frontal slices are M-matrices, we propose a tensor splitting iterative algorithm. Starting from $\yv_0\in\Delta$, we repeatedly find a new iterate $\yv_{c+1}$ by solving the following subproblems
\begin{equation}\label{TS0}
  \left(I\circ\one^{\circ(k-2)}\right)\yv_{c+1}^{k-1} = \alpha\bar{\Pt}\yv_c^{k-1}+\vv
\end{equation}
for $c=0,1,2,\dots$. This structural tensor equation has a closed-form solution.

\begin{Lem}\label{Lem:3-1}
  Let $\yv_c\in\Delta$. Equation \eqref{TS0} has a closed-form solution:
  \begin{equation*}
    \yv_{c+1} = \left(1+\alpha\one^T(\bar{\Pt}\yv_c^{k-1})\right)^{-\frac{k-2}{k-1}}(\alpha\bar{\Pt}\yv_c^{k-1}+\vv) \in\Delta,
  \end{equation*}
  which is exactly $\yv_{c+1}=\Phi(\yv_c)$, where $\Phi$ is defined in Lemma \ref{Lem:3-5}.
\end{Lem}
\begin{proof}
  It follows from a subproblem \eqref{TS0} that
  \begin{equation}\label{TS1}
    (\one^T\yv_{c+1})^{k-2}\yv_{c+1} = \alpha\bar{\Pt}\yv_c^{k-1}+\vv,
  \end{equation}
  which means $\yv_{c+1}$ is parallel to the right-hand-side vector $\alpha\bar{\Pt}\yv_c^{k-1}+\vv$. Hence, we set
  \begin{equation}
    \yv_{c+1} = \gamma(\alpha\bar{\Pt}\yv_c^{k-1}+\vv)
  \end{equation}
  with a positive undetermined parameter $\gamma$.

  To determine $\gamma$, we take inner products between $\one$ and both sides of \eqref{TS1} and thus get
  \begin{equation*}
    \one^T(\alpha\bar{\Pt}\yv_c^{k-1}+\vv) = (\one^T\yv_{c+1})^{k-1}
    = \gamma^{k-1}\left(\one^T(\alpha\bar{\Pt}\yv_c^{k-1}+\vv)\right)^{k-1},
  \end{equation*}
  which means
  \begin{equation}
    \gamma = \left(\one^T(\alpha\bar{\Pt}\yv_c^{k-1}+\vv)\right)^{-\frac{k-2}{k-1}}
    = \left(1+\alpha\one^T(\bar{\Pt}\yv_c^{k-1})\right)^{-\frac{k-2}{k-1}}.
  \end{equation}
  The proof is completed.
\end{proof}

\begin{algorithm}[t]
\caption{Tensor splitting algorithm for solving MLPPR.}\label{Alg:TSA}
\begin{algorithmic}[1]
  \State Choose $\yv_0\in\Re^n_+$ such that $\one^T\yv_0\le (1-\alpha)^{-\frac{1}{k-1}}$. Set $c\gets0$.
  \While{the sequence of iterates does not converge}
    \State Compute $\zv_c=\alpha\bar{\Pt}\yv_c^{k-1}+\vv$.
    \State Calculate $\yv_{c+1} = \left(\one^T\zv_c\right)^{-\frac{k-2}{k-1}}\zv_c.$
    \State $c \gets c+1$.
  \EndWhile
  % \STATE Normalize $\yv$ to obtain a stochastic vector $\xv=\yv/(\one^T\yv)$.
\end{algorithmic}
\end{algorithm}

Now, we present formally the tensor splitting algorithm in Algorithm \ref{Alg:TSA}, which is a fixed-point iteration \cite{LGL17}. On the global convergence of the tensor splitting algorithm, we have the following theorem.

\begin{Thm}\label{Thm:4-2}
  Let $\varsigma$ be defined as in Theorem \ref{Thm:3-A}.
  If $\varsigma<1$, then the tensor splitting algorithm generates a sequence of iterates $\{\yv_c\}$ which converges globally to a nonnegative solution $\yv_*$ of MLPPR with a linear rate, i.e.,
  \begin{equation*}
    \|\yv_{c}-\yv_*\|_1 \le \varsigma^c \|\yv_0-\yv_*\|_1.
  \end{equation*}
\end{Thm}
\begin{proof}
  Since $\Phi$ is a contractive map, we have
  \begin{equation*}
    \|\yv_{c+1}-\yv_*\|_1 = \|\Phi(\yv_c)-\Phi(\yv_*)\|_1 \le \varsigma\|\yv_c-\yv_*\|_1 \le \cdots
    \le \varsigma^{c+1} \|\yv_0-\yv_*\|_1.
  \end{equation*}
  The proof is completed.
\end{proof}

\section{Hypergraph partitioning}

To evaluate the performance of our MLPPR model and the associated tensor splitting algorithm, we do some numerical experiments on hypergraph partitioning. The tensor splitting algorithm terminates if
\begin{equation*}
  \frac{\|\yv_c-\yv_{c-1}\|_{\infty}}{\|\yv_c\|_1} \le 10^{-8} \quad\text{ or }\quad
  \frac{\|(\one^T\yv_c)^{k-2}\yv_c-\alpha\bar{\Pt}\yv_c^{k-1}-\vv\|_{\infty}}{\|\yv_c\|_1^{k-1}} \le 10^{-10}.
\end{equation*}
Theorem \ref{Thm:4-2} ensures that the tensor splitting algorithm for solving MLPPR converges if $(2k-3)\alpha(1-\alpha)^{-\frac{k-2}{k-1}}<1$. Not only that, for MLPPR with $\alpha\in[0,1)$, the tensor splitting algorithm always converges in our experiments.

\subsection{A toy example}\label{toy-hgraph}

% This subsection aims to quantify variations in the nonnegative solution $\yv$ of MLPPR with respect to perturbations in both the columnwise-substochastic tensor $\bar{\Pt}$ and the stochastic vector $\vv$.
% Recalling Theorem \ref{Thm:3-D}, we will identify which perturbation (in $\bar{\Pt}$ or $\vv$) has more impact on $\yv$.
Recalling Theorem \ref{Thm:3-D}, there is a theoretical bound of the variation in the MLPPR solution $\yv$ with respect to perturbations in both the columnwise-substochastic tensor $\bar{\Pt}$ and the stochastic vector $\vv$. Now, we aim to identify whether this theoretical bound is tight by a numerical way.

\begin{figure}[t]
  \centering
  \includegraphics[width=0.7\textwidth]{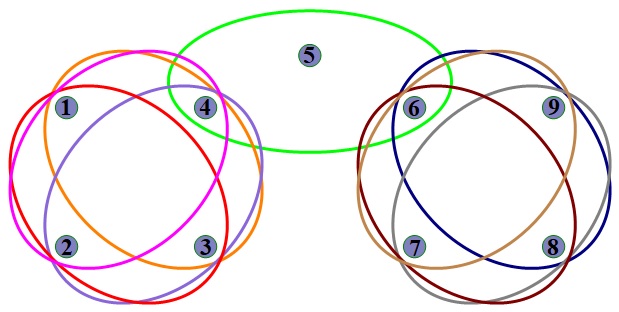}
  \caption{A toy example of a $3$-uniform hypergraph.}\label{ToyEx-HGraph}
\end{figure}

\begin{figure}[t]
  %\hspace{-60pt}
  \includegraphics[width=0.98\textwidth]{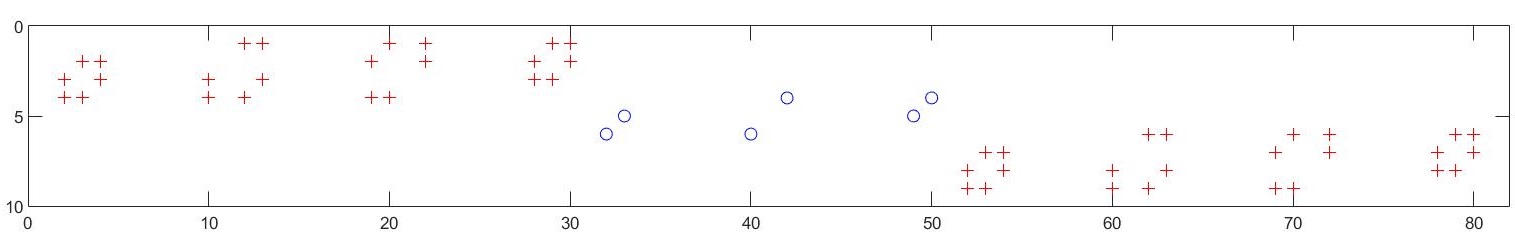}
  \caption{The mode-1 unfolding $\Rm(\bar{\Pt})\in\Re^{9\times 81}$ of the columnwise-substochastic tensor $\bar{\Pt}$ is illustrated, where a blue circle and a red plus stand for $1$ and $1/2$ respectively.}\label{spy-Pt}
\end{figure}

We consider an undirected $3$-uniform hypergraph $G=(V,E,\wv)$ illustrated in Fig. \ref{ToyEx-HGraph}, where $V=\{1,2,\dots,9\}$, $E=\left\{\{1,2,3\},\{1,2,4\},\{1,3,4\},\right.$ $\left.\{2,3,4\},\{4,5,6\},\{6,7,8\},\{6,7,9\},\{6,8,9\},\{7,8,9\}\right\}$, and $\wv=\one_9$. The adjacency tensor $\At$ and the columnwise-substochastic tensor $\bar{\Pt}$ of this hypergraph are produced by \eqref{AdjTen-undir} and \eqref{StoT}, respectively. Readers may refer to Fig. \ref{spy-Pt} for the sparsity of $\bar{\Pt}$. It is valuable to see that 51 out of 81 fibers of $\bar{\Pt}$ are dangling. To find a largest complete sub-hypergraph in $G$, we take a stochastic vector $\vv=(\frac{1}{2},\frac{1}{2},0,\dots,0)^T\in\Re^9_+$ and a probability $\alpha=\frac{1}{5}$. By solving MLPPR with these $\bar{\Pt}$, $\vv$, and $\alpha$, we obtain the nonnegative solution $\yv=(0.4796,0.4796,0.0527,0.0527,0.0000,\dots,0.0000)^T\in\Re^9$.

Next, to perturb a stochastic vector $\uv\in\Re^n$ ($\uv\ge\mathbf{0}$ and $\one^T\uv=1$) in a quantity $\sigma\in(0,1)$, we utilize an optimization approach for producing a random perturbation $\delta\uv\in\Re^n$ such that $\|\delta\uv\|_1=\sigma$, $\uv+\delta\uv\ge\mathbf{0}$, and $\one^T(\uv+\delta\uv)=1$.
First of all, a random vector $\bv\in\Re^n$ is generated by using the standard normal distribution $N(0,1)$. Then vector $\delta\uv$ takes the value from the projection of this random vector $\bv$ onto a compact set:
\begin{equation*}
  \min_{\delta\uv\in\Re^n}~ \|\delta\uv-\bv\|_2^2 \quad \mathrm{s.t.}~~\|\delta\uv\|_1=\sigma,~~
  \one^T\delta\uv = 0,~~\delta\uv \ge -\uv.
  %\begin{aligned}
  %  \min~ & \|\delta\uv-\bv\|_2^2 \\
  %  \mathrm{s.t.}~~ & \|\delta\uv\|_1=\sigma, \\
  %       & \one^T\delta\uv = 0, \\
  %       & \delta\uv \ge -\vv.
  %\end{aligned}
\end{equation*}
By introducing variables $\sv_+=\max(\delta\uv,\mathbf{0})\in\Re^n_+$ and $\sv_-=\sv_+-\delta\uv\in\Re^n_+$, the above projection problem could be represented as a convex quadratic programming
\begin{equation*}
  \begin{aligned}
    \min~ & \|\sv_+-\sv_- -\bv\|_2^2 \\
    \mathrm{s.t.}~~ & \one^T\sv_+ + \one^T\sv_- = \sigma, ~~ \one^T\sv_+ - \one^T\sv_- = 0, \\
         & \sv_+ - \sv_- \ge -\uv, ~~\sv_+ \ge\mathbf{0},~~  \sv_- \ge\mathbf{0},
  \end{aligned}
\end{equation*}
which could be solved efficiently. Once the optimal solution $(\sv_+,\sv_-)$ of this convex quadratic programming is available, it is straightforward to calculate $\delta\uv=\sv_+-\sv_-$.

\begin{figure}[t]
  \centering
  \includegraphics[width=0.8\textwidth]{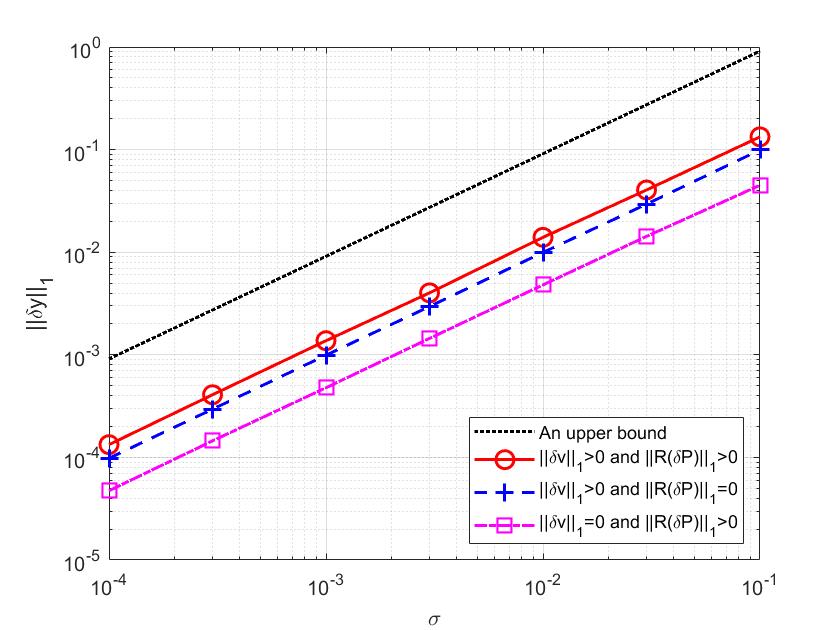}
  \caption{MLPPR solutions vary with various perturbations in the columnwise-substochastic tensor $\bar{\Pt}$ and a stochastic vector $\vv$.}\label{Perturb}
\end{figure}

Owing to the difference in coefficients of $\|\delta\vv\|_1$ and $\|\Rm(\delta\bar{\Pt})\|_1$ in Theorem \ref{Thm:3-D}, we set $\|\delta\vv\|_1=\sigma$ and $\|\Rm(\delta\bar{\Pt})\|_1=4\sigma$ to balance impacts from $\vv$ and $\bar{\Pt}$, where $\sigma$ is a perturbation quantity. For the columnwise-substochastic tensor $\bar{\Pt}$, only non-dangling fibers are perturbed. When the perturbation quantity $\sigma$ varies from $0.0001$ to $0.1$, we perform one hundred tests for each $\sigma$ and evaluate the mean of resulting changes $\|\delta\yv\|_1$ in the MLPPR solution. For comparison, we examine three cases of perturbations: (i) perturb both $\delta\vv$ and $\delta\bar{\Pt}$, (ii) only perturb $\delta\vv$, and (iii) only perturb $\delta\bar{\Pt}$. Numerical results on $\|\delta\yv\|_1$ versus $\sigma$ are illustrated in Figure \ref{Perturb}. The theoretical bound in Theorem \ref{Thm:3-D} is also printed here. It is easy to see from Figure \ref{Perturb} that the theoretical bound for $\|\delta\yv\|_1$ are almost tight up to some constant multiples. In addition, perturbation $\delta\bar{\Pt}$ has less impact on the MLPPR solution $\yv$ than perturbation $\delta\vv$.

\subsection{Subspace clustering}

Subspace clustering \cite{BP'13,CQZ'17} is a crucial problem in artificial intelligence. A typical problem is to extract subsets of collinear points from a set of points in $\Re^d$, where $d\ge2$. In this setting, classical pairwise approaches cannot work. % because any pair of points defines a straight line, and
Higher-order similarity of collinear points as well as tensor representions should be exploited.

At the beginning, to identify collinear points, there exists an approach \cite{GeoT'99} for computing a least-squares error of fitting a straight line with distinct points $\uv_1,\uv_2,\dots,\uv_k\in\Re^d$, where $k\ge3$.
% a least-squares method is adopted for finding the least-squares fitting line for given $s$ $(s\ge3)$ distinct points $\{\uv_1,\uv_2,\dots,\uv_s\}\subset\Re^d$.
% Let $\bar{\uv}:=\frac{1}{k}\sum_{i=1}^k\uv_i$ be the centre of points $\{\uv_1,\uv_2,\dots,\uv_s\}$.
% The best fitting line must goes through $\bar{\uv}$.
Assume that the straight line goes through $\uv_0$ with a unit direction $\vv$. The cost of least-squares fitting is
\begin{equation*}
  \Psi(\uv_0,\vv) := \sum_{i=1}^k \left\|(I-\vv\vv^T)(\uv_i-\uv_0)\right\|^2
    = \sum_{i=1}^k (\uv_i-\uv_0)^T(I-\vv\vv^T)(\uv_i-\uv_0).
\end{equation*}
To find an optimal $\uv_0$, we consider the partial derivative of the cost $\Psi$ with respect to $\uv_0$:
\begin{equation*}
  \frac{\partial \Psi(\uv_0,\vv)}{\partial \uv_0} = -2(I-\vv\vv^T) \sum_{i=1}^k (\uv_i-\uv_0)=0,
\end{equation*}
of which solutions satisfy $\sum_{i=1}^k (\uv_i-\uv_0)=\varrho\vv$ for $\varrho\in\Re$. For simplicity, we set $\varrho=0$ and get $\widehat{\uv}_0=\frac{1}{k}\sum_{i=1}^k\uv_i$.
% A straight line of the least-squares fitting goes through the centre $\bar{\uv}:=\frac{1}{k}\sum_{i=1}^k\uv_i$ of theses points.
Furthermore, let $U:=[\uv_1-\widehat{\uv}_0,\dots,\uv_k-\widehat{\uv}_0]\in\Re^{d\times k}$. The cost $\Psi$ could be rewritten as
\begin{equation*}
  \Psi(\widehat{\uv}_0,\vv) = \mathrm{trace}(U^TU) - \vv^TUU^T\vv.
\end{equation*}
% Hence we assume the equation of the line as follows
%\begin{equation}\label{straightline}
%  \vv^T(\uv-\bar{\uv})=0,
%\end{equation}
%where $\vv\in\Re^d$ is an undetermined unit vector. The sum of squares of distances from points $\uv_1,\uv_2,\dots,\uv_s$ to the line \eqref{straightline} is
%\begin{equation*}
%  \sum_{i=1}^s (\vv^T(\uv_i-\bar{\uv}))^2 = \vv^TUU^T\vv.
%\end{equation*}
To minimize the cost $\Psi$, an optimal (unit) direction $\widehat\vv$ has a closed-form solution: the eigenvector corresponding to the largest eigenvalue of $UU^T$.

\begin{figure}[t]
  \centering
  \begin{tabular}{cc}
    \includegraphics[width=.46\textwidth]{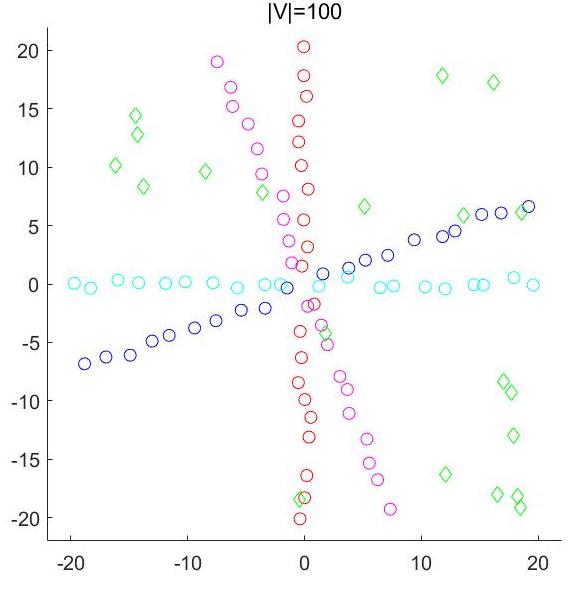} &
    \includegraphics[width=.46\textwidth]{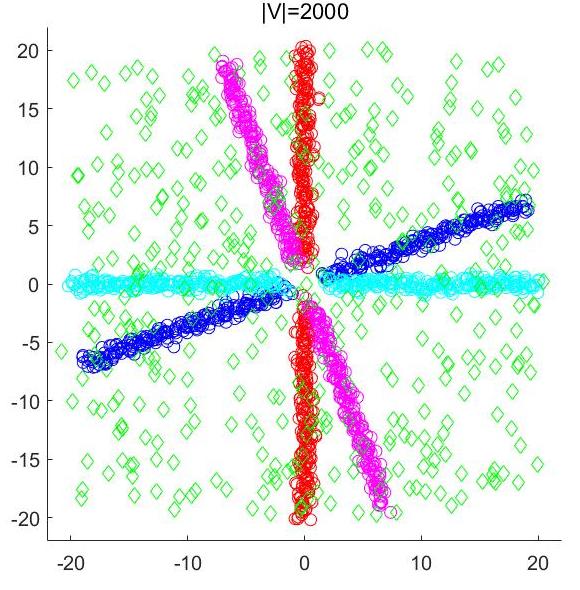}\\
    (a) & (b)
  \end{tabular}
  \caption{Illustrations of points for subspace clustering. Points of various colors (blue, cyan, magenta, and red) are located in different subspaces. Green diamonds are outliers.}\label{Example_SubSpace}
\end{figure}

Next, we consider clustering multiple subspaces. Fig. \ref{Example_SubSpace}(a) illustrates one hundred points, which constitute roughly four subsets of collinear points with blue, cyan, magenta, and red colors located in lines (subspaces) with slopes $\pi/9$, $0$, $-7\pi/18$, and $-\pi/2$, respectively. All points are corrupted by additive Gaussian noise $N(0,1/2)$. Each subset has twenty percent of points and the remainder are outliers that are marked as green diamonds. A larger example with two thousand points is displayed in Fig. \ref{Example_SubSpace}(b). When only positions of these points are known, how can we cluster these points such that points in each subset are nearly collinear?

As a first step, we construct an undirected $3$-uniform hypergraph $G=(V,E,\wv)$. Each point is regarded as a vertex in $V$. Hence, $\lvert V\rvert=n$. Basically, every three vertices may be connected as an edge to produce a complete $3$-uniform hypergraph. But it is too expensive to handle complete hypergraphs in practice. We prefer to use a random hypergraph which is generated by two steps economically. (i) To ensure connectedness, a complete graph is produced in advance. Then, for each edge of the complete graph, we add a different vertex which is randomly chosen. So there are $n(n-1)/2$ edges as candidates. (ii) For every candidate edge, we compute the least-squares error of fitting a straight line with positions of three points in a candidate edge. Then five percent candidate edges with smaller fitting errors are picked up to form $m:=\lvert E\rvert=[n(n-1)/40]$ edges of a random $3$-uniform hypergraph. All weights of edges of this hypergraph are ones.

\begin{algorithm}[t]
\caption{MLPPR-based hypergraph bipartition.}\label{Alg:CHmpPR}
\begin{algorithmic}[1]
  \State Form a $3$-uniform hypergraph $G$, generate its adjacency tensor $\At$ by \eqref{AdjTen-undir} and then a columnwise-substochastic tensor $\bar{\Pt}$ by \eqref{StoT}.
  \State Choose a probability $\alpha=0.99$ and a stochastic vector $\vv=\one/n$.
  \State Solve MLPPR \eqref{mlPPR} by Algorithm \ref{Alg:TSA} for a nonnegative solution $\yv_*$.
  \State Form $\widehat{A}=\bar{\Pt}\times_3\yv_*$ as an adjacency matrix of a latent directed graph.
  \State Compute $D=\mathrm{diag}(\one^T\widehat{A})$ and $\bar{P}=\widehat{A}D^{-1}$. Find $\mathbf{\pi}$ such that $\bar{P}\mathbf{\pi}=\mathbf{\pi}$.
  \State Let $\Pi=\mathrm{diag}(\mathbf{\pi})$. Find the second left eigenvector $\xv_*$ of $2P_{\mathrm{sym}}=\Pi\bar{P}^T\Pi^{-1}+\bar{P}$.
  \State Use $\xv_*$ to produce a bipartition of the vertex set of the given hypergraph $G$.
\end{algorithmic}
\end{algorithm}

Once a random hypergraph is ready, our approach, i.e., Algorithm \ref{Alg:CHmpPR} using MLPPR, follows the framework of Algorithm 1 in \cite{BGL'15}. Roughly speaking, there are two stages. (i) MLPPR-based Hypergraph analysis. Above all, we generate an adjacency tensor $\At$ and a columnwise-substochastic tensor $\bar{\Pt}$ of a random $3$-uniform hypergraph. By setting a probability $\alpha=0.99$ and a stochastic vector $\vv=\one/n$, an MLPPR problem is solved by the tensor splitting algorithm to obtain a nonnegative solution $\yv_*$.
(ii) Unsupervised clustering of a directed graph \cite{BGL'15,Ch'05}. Let $\widehat{A}:=\bar{\Pt}\times_3\yv_*\in\Re^{n\times n}$ be an adjacency matrix of a latent directed graph. It is straightforward to calculate a degree matrix $D=\mathrm{diag}(\one^T\widehat{A})$ and a columnwise-stochastic matrix $\bar{P}=\widehat{A}D^{-1}$. By a symmetrized approach in \cite{Ch'05}, we set $\Pi:=\mathrm{diag}(\mathbf{\pi})$, where $\mathbf{\pi}\in\Re^n_+$ is the stationary distribution of $\bar{P}$ such that $\bar{P}\mathbf{\pi}=\mathbf{\pi}$. Then, we define $A_{\mathrm{sym}}:=\frac{1}{2}(\Pi\bar{P}^T+\bar{P}\Pi)$, $D_{\mathrm{sym}}=\mathrm{diag}(A_{\mathrm{sym}}\one)=\Pi$, and $P_{\mathrm{sym}}=A_{\mathrm{sym}}D_{\mathrm{sym}}^{-1}=\frac{1}{2}(\Pi\bar{P}^T\Pi^{-1}+\bar{P})$. Finally, the second left eigenvector of $P_{\mathrm{sym}}$ is calculated for partitioning vertices of this hypergraph.

\begin{figure}
  \centering
  \includegraphics[width=0.8\textwidth]{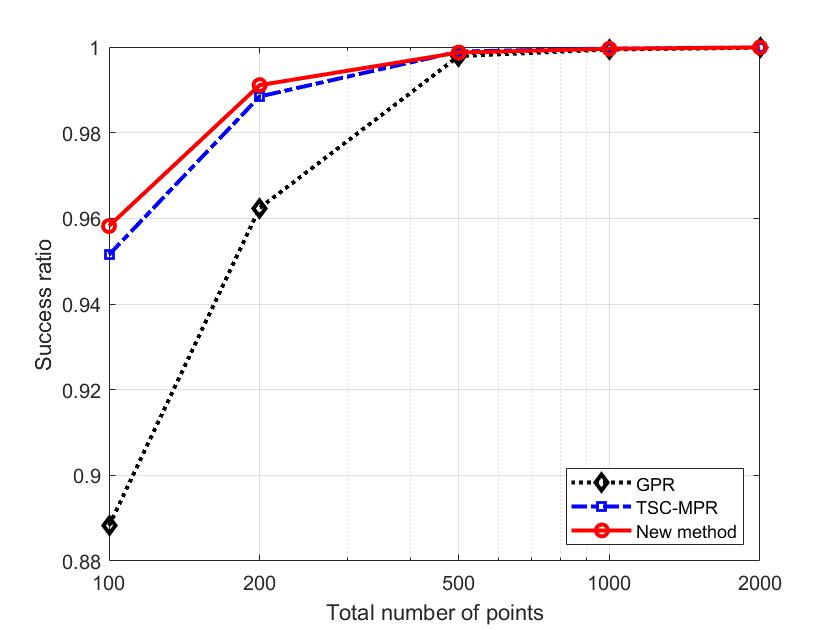}\\
  \caption{Success ratio of various methods}\label{Spaces:SussRatio}
\end{figure}

To produce a bipartition of the vertex set $V$ of a $3$-uniform hypergraph $G$, i.e., Step 7 of Algorithm \ref{Alg:CHmpPR}, we utilize the normalized cut of hypergraphs. For a vertex subset $\emptyset\subset S\subset V$, we define the volume of $S$:
\begin{equation*}
  \mathrm{vol}(S):=\sum_{i,j,k\atop i\in S}a_{ijk},
\end{equation*}
where $a_{ijk}$'s are elements of an adjacency tensor of $G$. Let $\bar{S}:=V\backslash S$. The cut of $S$ is defined as
\begin{equation*}
  \mathrm{cut}(S):=\sum_{i,j,k}a_{ijk}-\sum_{i,j,k\in S}a_{ijk}-\sum_{i,j,k\in \bar{S}}a_{ijk}.
\end{equation*}
Then, the normalized cut of $G$ is defined as follows
\begin{equation*}
  h = \min_{\emptyset\subset S\subset V} \mathrm{cut}(S)\left(\frac{1}{\mathrm{vol}(S)}+\frac{1}{\mathrm{vol}(\bar{S})}\right).
\end{equation*}
Indeed, it is NP-hard to evaluate normalized cuts of hypergraphs exactly. A heuristic and economic approach is widely used. We sort vertices in $V$ of $G$ according to the corresponding elements of the second left eigenvector $\xv_*$ of $P_{\mathrm{sym}}$ to obtain an ordered arrangement $\nu_1,\nu_2,\dots,\nu_n$ of vertices in $V$. To find an efficient cut, we only consider subsets $S$ that contains the leading vertices of this ordered arrangement. That is to say, we compute
\begin{equation}\label{h-ncut-0}
  h(i) = \mathrm{cut}(\{\nu_1,\dots,\nu_i\})\left(\frac{1}{\mathrm{vol}(\{\nu_1,\dots,\nu_i\})}
    +\frac{1}{\mathrm{vol}(\{\nu_{i+1},\dots,\nu_n\})}\right)
\end{equation}
for $i=1,\dots,n-1$ and then pick up the optimal index $i_*$ such that
\begin{equation}\label{h-ncut}
  i_* = \mathop{\arg\min}_{1\le i\le n-1} h(i).
\end{equation}
The minimization \eqref{h-ncut} is a heuristic approximation of the normalized cut of $G$. Then, a solution $i_*$ of \eqref{h-ncut} produces a bipartition of points corresponding to vertex subsets $\{\nu_1,\dots,\nu_{i_*}\}$ and $\{\nu_{i_*+1},\dots,\nu_n\}$ of $G$. Whereafter, since there are four subspaces in our experiments, we recursively run Algorithm \ref{Alg:CHmpPR} to repartition a sub-hypergraph with a larger vertex subset of $G$ until the set of points are clustered in four parts.

A significant difference between our method and Algorithm 1 in \cite{BGL'15} (named TSC-MPR) is that the former do not use any dangling information but the latter do. On one hand, for large-scale (sparse) random hypergraphs, there are a mount of dangling fibers contained in adjacency tensors of hypergraphs. By discarding dangling information, MLPPR is more efficient than MPR. On the other hand, although the stochastic vector $\yv_*/(\one^T\yv_*)$, where $\yv_*$ is generated by MLPPR, is the same as the solution $\xv_*$ of MPR, the matrix $\widehat{A}=\bar{\Pt}\times_3\yv_*$ in Step 4 of Algorithm \ref{Alg:CHmpPR} is different from $\Pt\times_3\xv_*$ used in TSC-MPR. Since the tensor $\Pt$ is a columnwise-stochastic tensor with dangling corrections, the matrix $\Pt\times_3\xv_*$ is also columnwise-stochastic. When we give up dangling corrections, the matrix $\widehat{A}=\bar{\Pt}\times_3\yv_*$ is only nonnegative but not necessarily columnwise-stochastic. In this sense, the matrix $\widehat{A}=\bar{\Pt}\times_3\yv_*$ is clean while the matrix $\Pt\times_3\xv_*$ may be corrupted by dangling corrections.

For solving subspace clustering problems, we run three methods: (i) The new method, i.e., Algorithm \ref{Alg:CHmpPR} using MLPPR. (ii) TSC-MPR: tensor spectral clustering by MPR in \cite{BGL'15}. (iii) GPR: PageRank for a graph approximation of a hypergraph. For subspace clustering problems with total numbers of points ranging from $100$ to $2000$, we count the success ratios of points being in correct clusters, where outliers could be in any clusters. For each case, we run 100 tests with randomly generated points. Numerical results illustrated in Fig. \ref{Spaces:SussRatio} means that the new method outperforms TSC-MPR and GPR. Owing to the skilful avoiding of dangling corrections, the new method is effective and efficient.

\subsection{Network analysis}

Higher-order graph information is important for analyzing network. For example, triangles (network structures involving three nodes) are fundamental to understand social network and community structure \cite{KW06,RELWL14,BGL16,BGL'17}. Particularly, we are interested in the directed 3-cycle (D3C), which captures a feedback loop in network.

\begin{table}
  \begin{center}
  \begin{minipage}{0.85\textwidth}
  \caption{Information of networks and partition sizes.}\label{Net:1}
  \begin{tabular}{lrr|ccc}
    \hline
    Network       & \#nodes & \#D3Cs    & GPR & TSC-MPR & New method  \\
    \hline
    email-EuAll   & 11,315  & 183,836   & 1,492 & 1,700 & 1,634 \\
    % soc-Epinionsl & 15,963  & 738,231   & & & \\
    wiki-Talk     & 52,411  & 5,138,613 & 21,893 & 21,240 & 21,238 \\
    \hline
  \end{tabular}
  \end{minipage}
  \end{center}
\end{table}

We consider two practical networks \cite{snapnets} in Table \ref{Net:1} from \path"https://snap.stanford.edu/data/". Originally, each network is a directed graph. For instance, each node in the ``email-EuAll'' network represents an email address. There exists a directed edge from nodes $i$ to $j$ if $i$ sent at least one message to $j$. Before higher-order network analysis, we filter the network as follows: (i) remove all vertices and all edges that do not participate in any D3Cs; (ii) take the largest strongly connected component of the remaining network. As a result, the filtered ``email-EuAll'' network contains $11,315$ nodes and $183,836$ D3Cs. A larger network named ``wiki-Talk'' contains $52,411$ nodes and $5,138,613$ D3Cs.

\begin{figure}
  \centering
  \begin{tabular}{cc}
    \includegraphics[width=0.46\textwidth]{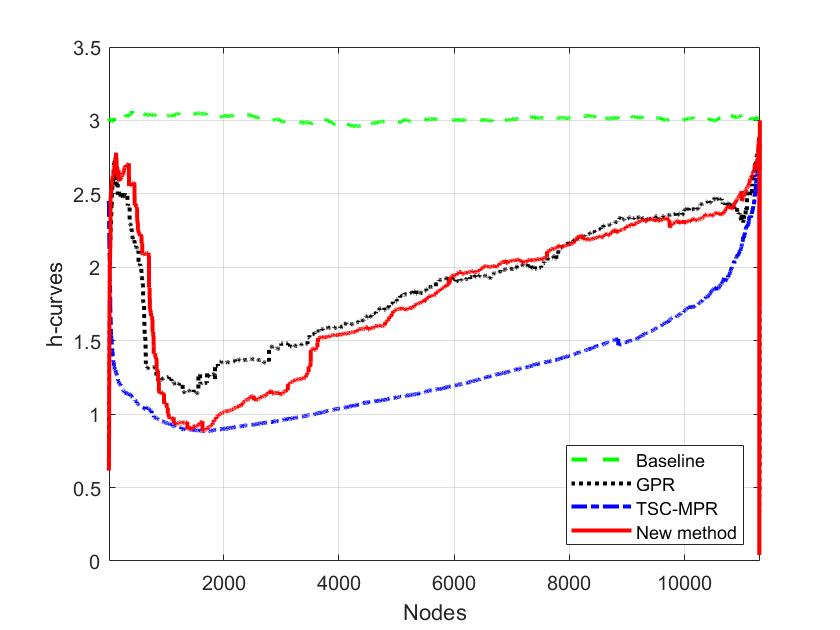} & \includegraphics[width=0.46\textwidth]{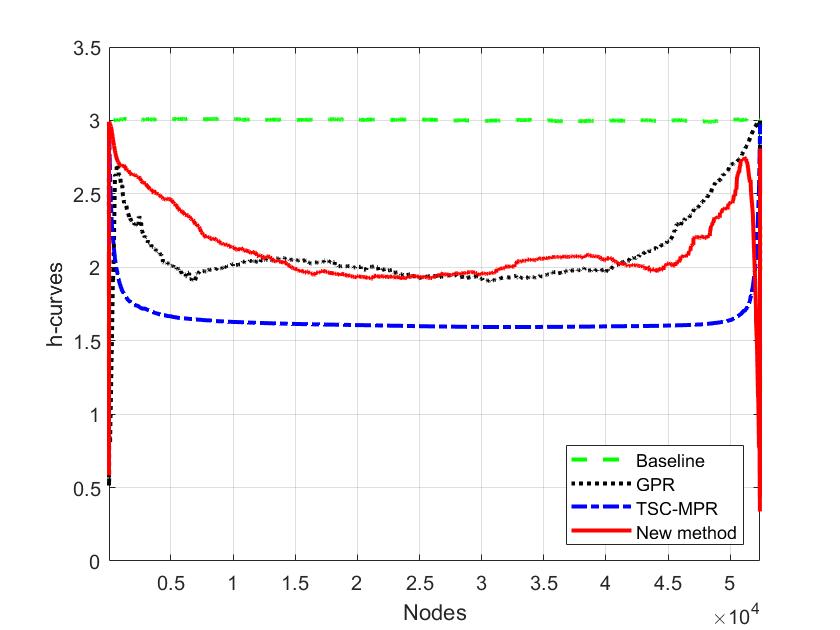} \\
    (a) ``email-EuAll'' network & (b) ``wiki-Talk'' network
  \end{tabular}
  \caption{$h$-curves of networks.}\label{NA:NCUT}
\end{figure}

\begin{figure}
  \centering
  \begin{tabular}{ccc}
    \includegraphics[width=0.3\textwidth]{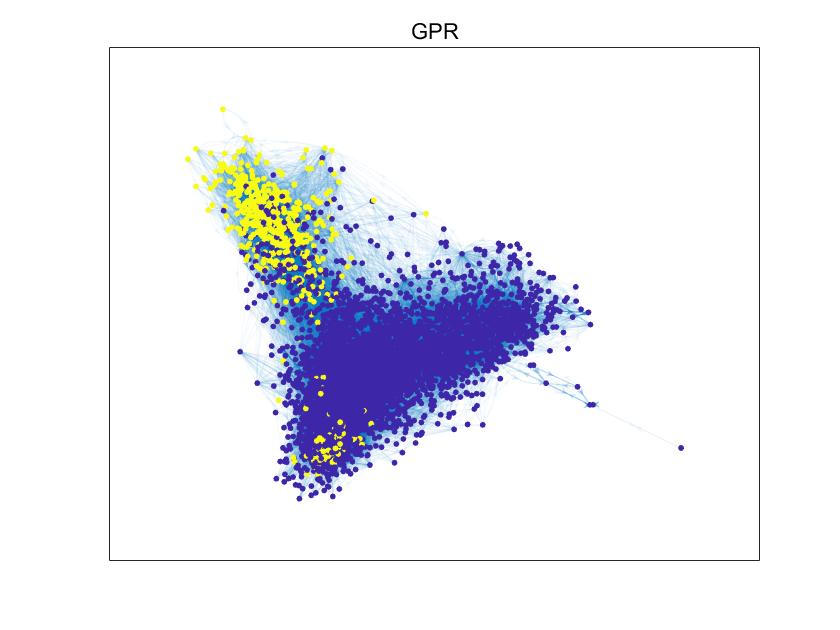} &
    \includegraphics[width=0.3\textwidth]{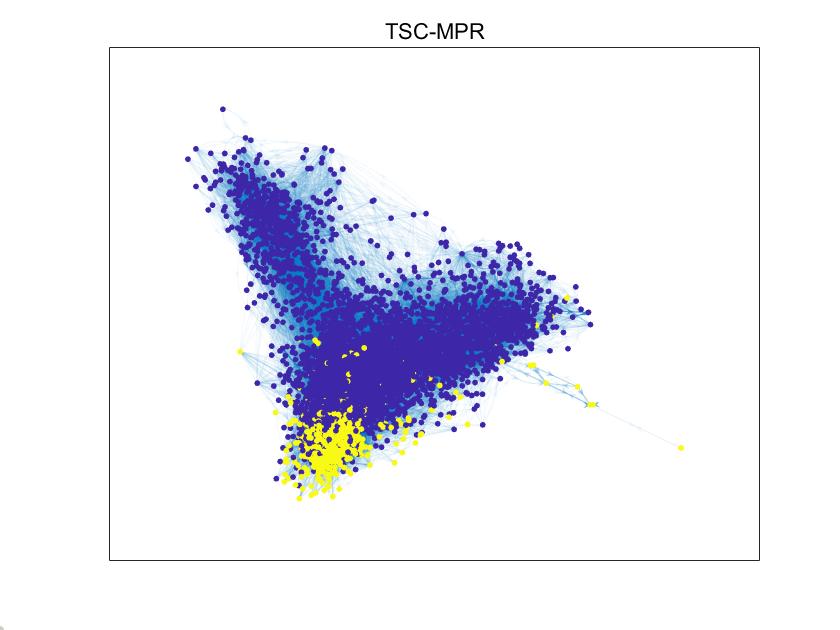} &
    \includegraphics[width=0.3\textwidth]{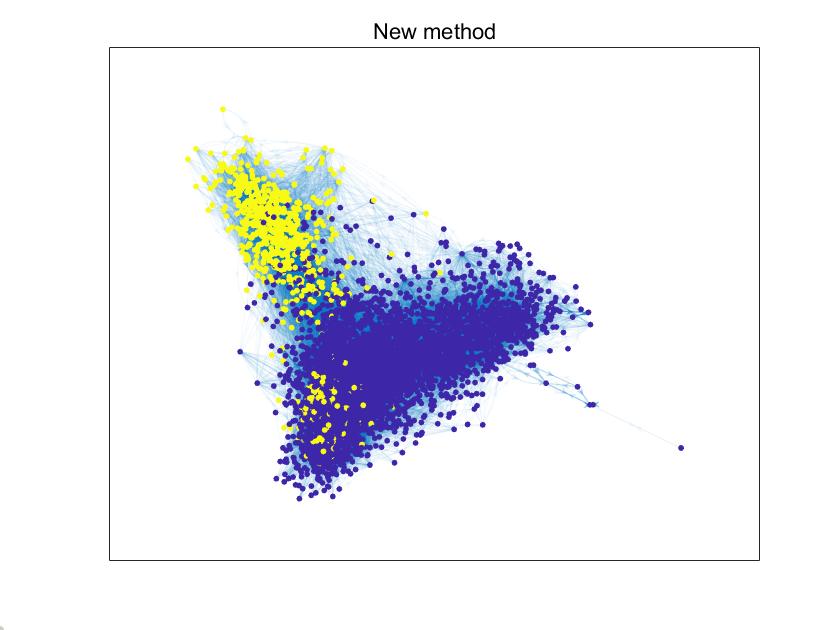} \\
    (a) & (b) & (c) \\
    \includegraphics[width=0.3\textwidth]{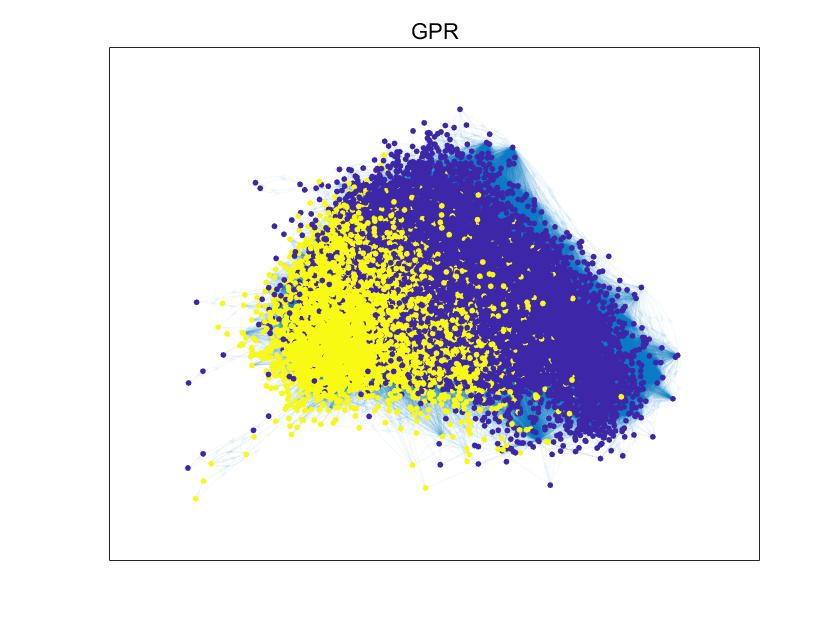} &
    \includegraphics[width=0.3\textwidth]{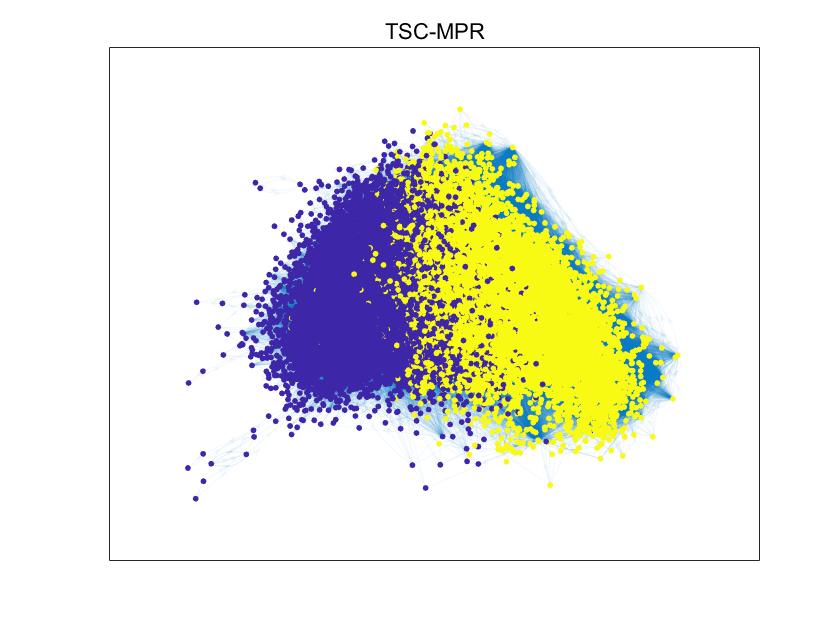} &
    \includegraphics[width=0.3\textwidth]{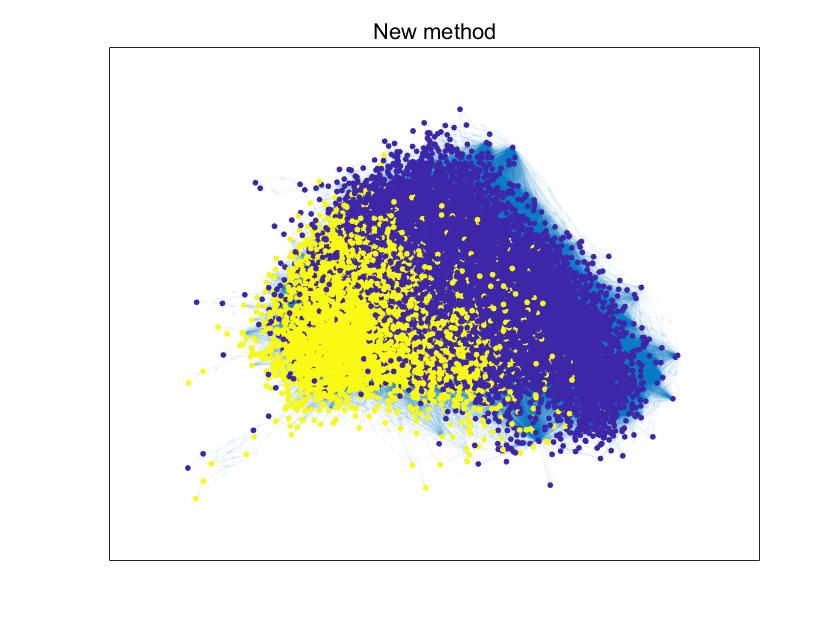} \\
    (d) & (e) & (f) \\
  \end{tabular}
  \caption{Illustration of partitioned results. Results of the ``email-EuAll'' network are displayed in (a--c) and results of the ``wiki-Talk'' network are displayed in (d--f). }\label{NA:Netw}
\end{figure}

Then it is straightforward to construct a $3$-uniform hypergraph $G$, where nodes and D3Cs of a network are vertices and edges of a $3$-uniform hypergaph, respectively. All weights are set to ones. To produce heuristic and efficient cuts, we test four algorithms for generating ordered arrangements of vertices of a hypergraph: (i) Baseline: a random order, (ii) GPR: PageRank for a graph approximation of a hypergraph, (iii) TSC-MPR: tensor spectral clustering by MPR in \cite{BGL'15}, and (iv) the new method, i.e., Algorithm \ref{Alg:CHmpPR} using MLPPR. By experiments, four algorithms generate four orders, of which associated $h$-curves defined by \eqref{h-ncut-0} are illustrated in Fig. \ref{NA:NCUT}. By picking up a minimum of $h$ by \eqref{h-ncut}, each  algorithm obtains a partition. The new method returns bipartitions with $\lvert S\rvert=1,634$ nodes for the ``email-EuAll'' network and $\lvert S\rvert=21,238$ nodes for the ``wiki-Talk'' network, as illustrated respectively in Fig. \ref{NA:Netw} (c) and (f). Partitions produced by TSC-MPR and GPR are also reported in Fig. \ref{NA:Netw} and associated partition sizes are listed in Table \ref{Net:1}. When compared with these existing methods, the new method is competitive.

\begin{table}[t]
  \begin{center}
  \begin{minipage}{0.5\textwidth}
  \caption{CPU times (seconds).}\label{NA:cputime}
  \centering
\begin{tabular}{c|ccc}
  \hline
  Network       & MPR-1  & MPR-2   & MLPPR \\
  \hline
  email-EuAll   &  1.882 &  0.073 & 0.064 \\
  wiki-Talk     & 58.530 &  2.257 & 0.913 \\
  \hline
\end{tabular}
  \end{minipage}
  \end{center}
\end{table}

Finally, we examine CPU times costed by solving MPR and MLPPR problems. MLPPR problems are solved by the tensor splitting algorithm and MPR problems are solved by the shifted fixed-point iteration \cite{BGL'15,GLY'15}. We should remind readers that CPU times costed by any existing algorithms for solving MPR and MLPPR are just a little part of the time for handling the whole process of network analysis. We recall the dangling correction \eqref{dePseudo} used by MPR:
\begin{equation*}
  \Pt := \bar{\Pt}+\vv\circ \underbrace{\left(\one^{\circ(k-1)}-\bar{\Pt}\bar{\times}_1 \one\right)}_{:=\mathcal{C}} \in\Re^{[k,n]}_+,
\end{equation*}
where $\mathcal{C}$
% MPR in \eqref{mlPR}, the columnwise-stochastic tensor $\Pt$ therein is a correction of the normalized adjacency tensor of the hypergraph $G$ that satisfies
%\begin{equation*}
%  \Rm(\Pt) = \Rm(\At)D^+ + \vv\cv^T,
%\end{equation*}
%where $D=\mathrm{diag}(\e^T\Rm(\At))$ is nonnegative, $\Rm(\At)D^+$ is only columnwise-substochastoc, $\cv^T=\one_{n^2}^T-\one^T\Rm(\At)D^+$
is a $(k-1)$th order indicator tensor of dangling fibers of $\bar{\Pt}$. On one hand, a straightforward approach (MPR-1) is to save and to manipulate $\mathcal{C}$ explicitly.
% We note that indicator tensors $\mathcal{C}$ may have more nonzero elements than columnwise-substochastic tensors $\bar{\Pt}$ of large scale sparse hypergraphs.
On the other hand, if $\mathcal{C}$ is dense but $\one^{\circ(k-1)}-\mathcal{C}$ is sparse, an alternative approach (MPR-2) utilizes $\one^{\circ(k-1)}-\mathcal{C}$ explicitly. For solving MPR/MLPPR problems raising from ``email-EuAll'' and ``wiki-Talk'' networks, CPU times (in seconds) of MPR-1, MPR-2, and MLPPR are reported in Table \ref{NA:cputime}. Obviously, MPR-2 is faster than MPR-1. This verifies that indicator tensors $\mathcal{C}$ are dense. Particularly, MLPPR is about 60 times faster than MPR-1 when dealing with ``wiki-Talk'' network. According to these numerical experiments, the proposed MLPPR is powerful and efficient.

\section{Conclusions}

As an interesting extension of the classical PageRank model for graph partitioning, we proposed the MLPPR model for partitioning uniform hypergraphs. By using Laplacian tensors of uniform hypergraphs without any dangling corrections, MLPPR is more suitable for processing large-scale hypergraphs. The existence and uniqueness of nonnegative solutions of MLPPR were analyzed. An error bound of the MLPPR solution was analyzed when the Laplacian tensor and a stochastic vector are slightly perturbed. By exploiting the structural Laplacian tensors, we designed a tensor splitting algorithm (a fixed-point iterations) for MLPPR. Numerical experiments on subspace clustering and network analysis verified the powerfulness and effectiveness of our method.
%However, it is well known that the solution of a system of higher-order polynomials is usually unstable in the viewpoint of numerical computations. As a special system of polynomials, MLPPR can not get rid of this curse.

% Many optimization algorithms could be used to solve the least squares minimization with nonnegative constraints \eqref{mPR1} even though the teleportation parameter gets close to one. Based on {\L}ojasiewicz property of the tensor optimization problem, we customized and analyzed the gradient projection algorithm. On one hand, since MLPPR always has a solution, the global optimal value of the constrained least squares problem \eqref{mPR1} is zero. On the other hand, the gradient projection algorithm suffers from local critical points of which the objective function value are nonzero. We show that this gap could be filled in theory by running the gradient projection algorithm from finite many starting points and picking up the best solution. This interesting result is a natural corollary of the approximate form of Caratheodory's theorem.

%Declarations
%Conflict of interest The author states that he has no %financial interest in the subject matter of this research.
%There are no competing interests to declare.
%Ethical Standards Acceptance of the journal’s ethical %standards and agreement to follow them: The author
%accepts the journal’s ethical standards

\noindent \textbf{Funding} This work was funded by the National Natural Science Foundation of China grant 12171168, 12071159, 11901118, 11771405, 62073087 and U1811464.
\vspace{0.25cm}

\noindent \textbf{Data Availability} The data that support the findings of this study are available from the corresponding author upon reasonable request.

\section*{Declarations}

\textbf{Conflict of interest} The authors have no relevant financial or non-financial interests to disclose.
%Acknowledgments are not compulsory. Where included they should be brief. Grant or contribution numbers may be acknowledged.
%
%Please refer to Journal-level guidance for any specific requirements.


% \bibliography{mlppr_refs}% common bib file

\begin{thebibliography}{37}
\ifx \bisbn   \undefined \def \bisbn  #1{ISBN #1}\fi
\ifx \binits  \undefined \def \binits#1{#1}\fi
\ifx \bauthor  \undefined \def \bauthor#1{#1}\fi
\ifx \batitle  \undefined \def \batitle#1{#1}\fi
\ifx \bjtitle  \undefined \def \bjtitle#1{#1}\fi
\ifx \bvolume  \undefined \def \bvolume#1{\textbf{#1}}\fi
\ifx \byear  \undefined \def \byear#1{#1}\fi
\ifx \bissue  \undefined \def \bissue#1{#1}\fi
\ifx \bfpage  \undefined \def \bfpage#1{#1}\fi
\ifx \blpage  \undefined \def \blpage #1{#1}\fi
\ifx \burl  \undefined \def \burl#1{\textsf{#1}}\fi
\ifx \doiurl  \undefined \def \doiurl#1{\url{https://doi.org/#1}}\fi
\ifx \betal  \undefined \def \betal{\textit{et al.}}\fi
\ifx \binstitute  \undefined \def \binstitute#1{#1}\fi
\ifx \binstitutionaled  \undefined \def \binstitutionaled#1{#1}\fi
\ifx \bctitle  \undefined \def \bctitle#1{#1}\fi
\ifx \beditor  \undefined \def \beditor#1{#1}\fi
\ifx \bpublisher  \undefined \def \bpublisher#1{#1}\fi
\ifx \bbtitle  \undefined \def \bbtitle#1{#1}\fi
\ifx \bedition  \undefined \def \bedition#1{#1}\fi
\ifx \bseriesno  \undefined \def \bseriesno#1{#1}\fi
\ifx \blocation  \undefined \def \blocation#1{#1}\fi
\ifx \bsertitle  \undefined \def \bsertitle#1{#1}\fi
\ifx \bsnm \undefined \def \bsnm#1{#1}\fi
\ifx \bsuffix \undefined \def \bsuffix#1{#1}\fi
\ifx \bparticle \undefined \def \bparticle#1{#1}\fi
\ifx \barticle \undefined \def \barticle#1{#1}\fi
\bibcommenthead
\ifx \bconfdate \undefined \def \bconfdate #1{#1}\fi
\ifx \botherref \undefined \def \botherref #1{#1}\fi
\ifx \url \undefined \def \url#1{\textsf{#1}}\fi
\ifx \bchapter \undefined \def \bchapter#1{#1}\fi
\ifx \bbook \undefined \def \bbook#1{#1}\fi
\ifx \bcomment \undefined \def \bcomment#1{#1}\fi
\ifx \oauthor \undefined \def \oauthor#1{#1}\fi
\ifx \citeauthoryear \undefined \def \citeauthoryear#1{#1}\fi
\ifx \endbibitem  \undefined \def \endbibitem {}\fi
\ifx \bconflocation  \undefined \def \bconflocation#1{#1}\fi
\ifx \arxivurl  \undefined \def \arxivurl#1{\textsf{#1}}\fi
\csname PreBibitemsHook\endcsname

%%% 1
\bibitem{BP'98}
\begin{barticle}
\bauthor{\bsnm{Brin}, \binits{S.}},
\bauthor{\bsnm{Page}, \binits{L.}}:
\batitle{The anatomy of a large-scale hypertextual {W}eb search engine}.
\bjtitle{Computer Networks and ISDN Systems}
\bvolume{30}(\bissue{1}),
\bfpage{107}--\blpage{117}
(\byear{1998}).
\doiurl{10.1016/S0169-7552(98)00110-X}
\end{barticle}
\endbibitem

%%% 2
\bibitem{BL'06}
\begin{barticle}
\bauthor{\bsnm{Bryan}, \binits{K.}},
\bauthor{\bsnm{Leise}, \binits{T.}}:
\batitle{The $\$25,000,000,000$ eigenvector: The linear algebra behind
  {G}oogle}.
\bjtitle{SIAM Review}
\bvolume{48}(\bissue{3}),
\bfpage{569}--\blpage{581}
(\byear{2006}).
\doiurl{10.1137/050623280}
\end{barticle}
\endbibitem

%%% 3
\bibitem{Gl'15}
\begin{barticle}
\bauthor{\bsnm{Gleich}, \binits{D.F.}}:
\batitle{Page{R}ank beyond the {W}eb}.
\bjtitle{SIAM Review}
\bvolume{57}(\bissue{3}),
\bfpage{321}--\blpage{363}
(\byear{2015}).
\doiurl{10.1137/140976649}
\end{barticle}
\endbibitem

%%% 4
\bibitem{LM'06}
\begin{bbook}
\bauthor{\bsnm{Langville}, \binits{A.N.}},
\bauthor{\bsnm{Meyer}, \binits{C.D.}}:
\bbtitle{Google's {P}age{R}ank and Beyond: The Science of Search Engine
  Rankings}.
\bpublisher{Princeton university press},
\blocation{Princeton}
(\byear{2006})
\end{bbook}
\endbibitem

%%% 5
\bibitem{BGL16}
\begin{barticle}
\bauthor{\bsnm{Benson}, \binits{A.R.}},
\bauthor{\bsnm{Gleich}, \binits{D.F.}},
\bauthor{\bsnm{Leskovec}, \binits{J.}}:
\batitle{Higher-order organization of complex networks}.
\bjtitle{Science}
\bvolume{353}(\bissue{6295}),
\bfpage{163}--\blpage{166}
(\byear{2016}).
\doiurl{10.1126/science.aad9029}
\end{barticle}
\endbibitem

%%% 6
\bibitem{BGL'17}
\begin{barticle}
\bauthor{\bsnm{Benson}, \binits{A.R.}},
\bauthor{\bsnm{Gleich}, \binits{D.F.}},
\bauthor{\bsnm{Lim}, \binits{L.-H.}}:
\batitle{The spacey random walk: A stochastic process for higher-order data}.
\bjtitle{SIAM Review}
\bvolume{59}(\bissue{2}),
\bfpage{321}--\blpage{345}
(\byear{2017}).
\doiurl{10.1137/16M1074023}
\end{barticle}
\endbibitem

%%% 7
\bibitem{KW06}
\begin{barticle}
\bauthor{\bsnm{Kossinets}, \binits{G.}},
\bauthor{\bsnm{Watts}, \binits{D.J.}}:
\batitle{Empirical analysis of an evolving social network}.
\bjtitle{Science}
\bvolume{311}(\bissue{5757}),
\bfpage{88}--\blpage{90}
(\byear{2006}).
\doiurl{10.1126/science.1116869}
\end{barticle}
\endbibitem

%%% 8
\bibitem{RELWL14}
\begin{barticle}
\bauthor{\bsnm{Rosvall}, \binits{M.}},
\bauthor{\bsnm{Esquivel}, \binits{A.V.}},
\bauthor{\bsnm{Lancichinetti}, \binits{A.}},
\bauthor{\bsnm{West}, \binits{J.D.}},
\bauthor{\bsnm{Lambiotte}, \binits{R.}}:
\batitle{Memory in network flows and its effects on spreading dynamics and
  community detection}.
\bjtitle{Nature Communications}
\bvolume{5},
\bfpage{4630}
(\byear{2014}).
\doiurl{10.1038/ncomms5630}
\end{barticle}
\endbibitem

%%% 9
\bibitem{IS08}
\begin{barticle}
\bauthor{\bsnm{Ipsen}, \binits{I.C.F.}},
\bauthor{\bsnm{Selee}, \binits{T.M.}}:
\batitle{{PageRank} computation, with special attention to dangling nodes}.
\bjtitle{SIAM Journal on Matrix Analysis and Applications}
\bvolume{29}(\bissue{4}),
\bfpage{1281}--\blpage{1296}
(\byear{2008}).
\doiurl{10.1137/060664331}
\end{barticle}
\endbibitem

%%% 10
\bibitem{F15}
\begin{barticle}
\bauthor{\bsnm{Francesco}, \binits{T.}}:
\batitle{A note on certain ergodicity coefficients}.
\bjtitle{Special Matrices}
\bvolume{3}(\bissue{1}),
\bfpage{175}--\blpage{185}
(\byear{2015}).
\doiurl{10.1515/spma-2015-0016}
\end{barticle}
\endbibitem

%%% 11
\bibitem{LM06}
\begin{barticle}
\bauthor{\bsnm{Langville}, \binits{A.N.}},
\bauthor{\bsnm{Meyer}, \binits{C.D.}}:
\batitle{A reordering for the pagerank problem}.
\bjtitle{SIAM Journal on Scientific Computing}
\bvolume{27}(\bissue{6}),
\bfpage{2112}--\blpage{2120}
(\byear{2006}).
\doiurl{10.1137/040607551}
\end{barticle}
\endbibitem

%%% 12
\bibitem{Ch'05}
\begin{barticle}
\bauthor{\bsnm{Chung}, \binits{F.}}:
\batitle{Laplacians and the {C}heeger inequality for directed graphs}.
\bjtitle{Annals of Combinatorics}
\bvolume{9}(\bissue{1}),
\bfpage{1}--\blpage{19}
(\byear{2005}).
\doiurl{10.1007/s00026-005-0237-z}
\end{barticle}
\endbibitem

%%% 13
\bibitem{ACL'07}
\begin{barticle}
\bauthor{\bsnm{Andersen}, \binits{R.}},
\bauthor{\bsnm{Chung}, \binits{F.}},
\bauthor{\bsnm{Lang}, \binits{K.}}:
\batitle{Using {P}age{R}ank to locally partition a graph}.
\bjtitle{Internet Mathematics}
\bvolume{4}(\bissue{1}),
\bfpage{35}--\blpage{64}
(\byear{2007})
\end{barticle}
\endbibitem

%%% 14
\bibitem{ACL'08}
\begin{barticle}
\bauthor{\bsnm{Andersen}, \binits{R.}},
\bauthor{\bsnm{Chung}, \binits{F.}},
\bauthor{\bsnm{Lang}, \binits{K.}}:
\batitle{Local partitioning for directed graphs using {P}age{R}ank}.
\bjtitle{Internet Mathematics}
\bvolume{5}(\bissue{1-2}),
\bfpage{3}--\blpage{22}
(\byear{2008}).
\doiurl{im/1259158595}
\end{barticle}
\endbibitem

%%% 15
\bibitem{LN'14}
\begin{barticle}
\bauthor{\bsnm{Li}, \binits{W.}},
\bauthor{\bsnm{Ng}, \binits{M.K.}}:
\batitle{On the limiting probability distribution of a transition probability
  tensor}.
\bjtitle{Linear and Multilinear Algebra}
\bvolume{62}(\bissue{3}),
\bfpage{362}--\blpage{385}
(\byear{2014}).
\doiurl{10.1080/03081087.2013.777436}
\end{barticle}
\endbibitem

%%% 16
\bibitem{GLY'15}
\begin{barticle}
\bauthor{\bsnm{Gleich}, \binits{D.F.}},
\bauthor{\bsnm{Lim}, \binits{L.-H.}},
\bauthor{\bsnm{Yu}, \binits{Y.}}:
\batitle{Multilinear {P}age{R}ank}.
\bjtitle{SIAM Journal on Matrix Analysis and Applications}
\bvolume{36}(\bissue{4}),
\bfpage{1507}--\blpage{1541}
(\byear{2015}).
\doiurl{10.1137/140985160}
\end{barticle}
\endbibitem

%%% 17
\bibitem{CZ'13}
\begin{barticle}
\bauthor{\bsnm{Chang}, \binits{K.C.}},
\bauthor{\bsnm{Zhang}, \binits{T.}}:
\batitle{On the uniqueness and non-uniqueness of the positive {Z}-eigenvector
  for transition probability tensors}.
\bjtitle{Journal of Mathematical Analysis and Applications}
\bvolume{408}(\bissue{2}),
\bfpage{525}--\blpage{540}
(\byear{2013}).
\doiurl{10.1016/j.jmaa.2013.04.019}
\end{barticle}
\endbibitem

%%% 18
\bibitem{FT'20}
\begin{barticle}
\bauthor{\bsnm{Fasino}, \binits{D.}},
\bauthor{\bsnm{Tudisco}, \binits{F.}}:
\batitle{Ergodicity coefficients for higher-order stochastic processes}.
\bjtitle{SIAM Journal on Mathematics of Data Science}
\bvolume{2}(\bissue{3}),
\bfpage{740}--\blpage{769}
(\byear{2020}).
\doiurl{10.1137/19M1285214}
\end{barticle}
\endbibitem

%%% 19
\bibitem{LLNV'17}
\begin{barticle}
\bauthor{\bsnm{Li}, \binits{W.}},
\bauthor{\bsnm{Liu}, \binits{D.}},
\bauthor{\bsnm{Ng}, \binits{M.K.}},
\bauthor{\bsnm{Vong}, \binits{S.-W.}}:
\batitle{The uniqueness of multilinear {P}age{R}ank vectors}.
\bjtitle{Numerical Linear Algebra with Applications}
\bvolume{24}(\bissue{6}),
\bfpage{2107}--\blpage{112}
(\byear{2017}).
\doiurl{10.1002/nla.2107}
\end{barticle}
\endbibitem

%%% 20
\bibitem{LLVX'20}
\begin{barticle}
\bauthor{\bsnm{Li}, \binits{W.}},
\bauthor{\bsnm{Liu}, \binits{D.}},
\bauthor{\bsnm{Vong}, \binits{S.-W.}},
\bauthor{\bsnm{Xiao}, \binits{M.}}:
\batitle{Multilinear {P}age{R}ank: Uniqueness, error bound and perturbation
  analysis}.
\bjtitle{Applied Numerical Mathematics}
\bvolume{156},
\bfpage{584}--\blpage{607}
(\byear{2020}).
\doiurl{10.1016/j.apnum.2020.05.022}
\end{barticle}
\endbibitem

%%% 21
\bibitem{LCN'13}
\begin{barticle}
\bauthor{\bsnm{Li}, \binits{W.}},
\bauthor{\bsnm{Cui}, \binits{L.-B.}},
\bauthor{\bsnm{Ng}, \binits{M.K.}}:
\batitle{The perturbation bound for the {P}erron vector of a transition
  probability tensor}.
\bjtitle{Numerical Linear Algebra with Applications}
\bvolume{20}(\bissue{6}),
\bfpage{985}--\blpage{1000}
(\byear{2013}).
\doiurl{10.1002/nla.1886}
\end{barticle}
\endbibitem

%%% 22
\bibitem{B'19}
\begin{barticle}
\bauthor{\bsnm{Benson}, \binits{A.R.}}:
\batitle{Three hypergraph eigenvector centralities}.
\bjtitle{SIAM Journal on Mathematics of Data Science}
\bvolume{1}(\bissue{2}),
\bfpage{293}--\blpage{312}
(\byear{2019}).
\doiurl{10.1137/18M1203031}
\end{barticle}
\endbibitem

%%% 23
\bibitem{BGL'15}
\begin{bchapter}
\bauthor{\bsnm{Benson}, \binits{A.R.}},
\bauthor{\bsnm{Gleich}, \binits{D.F.}},
\bauthor{\bsnm{Leskovec}, \binits{J.}}:
\bctitle{Tensor Spectral Clustering for Partitioning Higher-order Network
  Structures}.
In: \bbtitle{Proc SIAM Int Conf Data Min.},
pp. \bfpage{118}--\blpage{126}
(\byear{2015})
\end{bchapter}
\endbibitem

%%% 24
\bibitem{MP'18}
\begin{barticle}
\bauthor{\bsnm{Meini}, \binits{B.}},
\bauthor{\bsnm{Poloni}, \binits{F.}}:
\batitle{Perron-based algorithms for the multilinear {P}age{R}ank}.
\bjtitle{Numerical Linear Algebra with Applications}
\bvolume{25}(\bissue{6}),
\bfpage{2177}--\blpage{115}
(\byear{2018}).
\doiurl{10.1002/nla.2177}
\end{barticle}
\endbibitem

%%% 25
\bibitem{HW'21}
\begin{barticle}
\bauthor{\bsnm{Huang}, \binits{J.}},
\bauthor{\bsnm{Wu}, \binits{G.}}:
\batitle{Convergence of the fixed-point iteration for multilinear
  {P}age{R}ank}.
\bjtitle{Numerical Linear Algebra with Applications}
\bvolume{28}(\bissue{5}),
\bfpage{2379}
(\byear{2021}).
\doiurl{10.1002/nla.2379}
\end{barticle}
\endbibitem

%%% 26
\bibitem{LLV'19}
\begin{barticle}
\bauthor{\bsnm{Liu}, \binits{D.}},
\bauthor{\bsnm{Li}, \binits{W.}},
\bauthor{\bsnm{Vong}, \binits{S.-W.}}:
\batitle{Relaxation methods for solving the tensor equation arising from the
  higher-order {M}arkov chains}.
\bjtitle{Numerical Linear Algebra with Applications}
\bvolume{26}(\bissue{5}),
\bfpage{2260}
(\byear{2019}).
\doiurl{10.1002/nla.2260}
\end{barticle}
\endbibitem

%%% 27
\bibitem{CRT'20}
\begin{barticle}
\bauthor{\bsnm{Cipolla}, \binits{S.}},
\bauthor{\bsnm{Redivo-Zaglia}, \binits{M.}},
\bauthor{\bsnm{Tudisco}, \binits{F.}}:
\batitle{Extrapolation methods for fixed-point multilinear {P}age{R}ank
  computations}.
\bjtitle{Numerical Linear Algebra with Applications}
\bvolume{27}(\bissue{2}),
\bfpage{2280}
(\byear{2020}).
\doiurl{10.1002/nla.2280}
\end{barticle}
\endbibitem

%%% 28
\bibitem{YCO'21}
\begin{barticle}
\bauthor{\bsnm{Yuan}, \binits{A.}},
\bauthor{\bsnm{Calder}, \binits{J.}},
\bauthor{\bsnm{Osting}, \binits{B.}}:
\batitle{A continuum limit for the {P}age{R}ank algorithm}.
\bjtitle{European Journal of Applied Mathematics}
\bvolume{33}(\bissue{3}),
\bfpage{472}--\blpage{504}
(\byear{2022}).
\doiurl{10.1017/S0956792521000097}
\end{barticle}
\endbibitem

%%% 29
\bibitem{BP'13}
\begin{barticle}
\bauthor{\bsnm{Bul\`{o}}, \binits{S.R.}},
\bauthor{\bsnm{Pelillo}, \binits{M.}}:
\batitle{A game-theoretic approach to hypergraph clustering}.
\bjtitle{IEEE Transactions on Pattern Analysis and Machine Intelligence}
\bvolume{35}(\bissue{6}),
\bfpage{1312}--\blpage{1327}
(\byear{2013}).
\doiurl{10.1109/TPAMI.2012.226}
\end{barticle}
\endbibitem

%%% 30
\bibitem{CD'12}
\begin{barticle}
\bauthor{\bsnm{Cooper}, \binits{J.}},
\bauthor{\bsnm{Dutle}, \binits{A.}}:
\batitle{Spectra of uniform hypergraphs}.
\bjtitle{Linear Algebra and its Applications}
\bvolume{436}(\bissue{9}),
\bfpage{3268}--\blpage{3292}
(\byear{2012}).
\doiurl{10.1016/j.laa.2011.11.018}
\end{barticle}
\endbibitem

%%% 31
\bibitem{QL17book}
\begin{bbook}
\bauthor{\bsnm{Qi}, \binits{L.}},
\bauthor{\bsnm{Luo}, \binits{Z.}}:
\bbtitle{Tensor Analysis: Spectral Theory and Special Tensors}.
\bpublisher{SIAM},
\blocation{Philadelphia}
(\byear{2017})
\end{bbook}
\endbibitem

%%% 32
\bibitem{GCH22}
\begin{barticle}
\bauthor{\bsnm{Gao}, \binits{G.}},
\bauthor{\bsnm{Chang}, \binits{A.}},
\bauthor{\bsnm{Hou}, \binits{Y.}}:
\batitle{Spectral radius on linear $r$-graphs without expanded $k_{r+1}$}.
\bjtitle{SIAM Journal on Discrete Mathematics}
\bvolume{36}(\bissue{2}),
\bfpage{1000}--\blpage{1011}
(\byear{2022}).
\doiurl{10.1137/21M1404740}
\end{barticle}
\endbibitem

%%% 33
\bibitem{LN15}
\begin{barticle}
\bauthor{\bsnm{Li}, \binits{W.}},
\bauthor{\bsnm{Ng}, \binits{M.K.}}:
\batitle{Some bounds for the spectral radius of nonnegative tensors}.
\bjtitle{Numerische Mathematik}
\bvolume{130},
\bfpage{315}--\blpage{335}
(\byear{2015}).
\doiurl{10.1007/s00211-014-0666-5}
\end{barticle}
\endbibitem

%%% 34
\bibitem{LGL17}
\begin{barticle}
\bauthor{\bsnm{Liu}, \binits{C.-S.}},
\bauthor{\bsnm{Guo}, \binits{C.-H.}},
\bauthor{\bsnm{Lin}, \binits{W.-W.}}:
\batitle{Newton--{N}oda iteration for finding the {P}erron pair of a weakly
  irreducible nonnegative tensor}.
\bjtitle{Numerische Mathematik}
\bvolume{137},
\bfpage{63}--\blpage{90}
(\byear{2017}).
\doiurl{10.1007/s00211-017-0869-7}
\end{barticle}
\endbibitem

%%% 35
\bibitem{CQZ'17}
\begin{barticle}
\bauthor{\bsnm{Chen}, \binits{Y.}},
\bauthor{\bsnm{Qi}, \binits{L.}},
\bauthor{\bsnm{Zhang}, \binits{X.}}:
\batitle{The {F}iedler vector of a {L}aplacian tensor for hypergraph
  partitioning}.
\bjtitle{SIAM Journal on Scientific Computing}
\bvolume{39}(\bissue{6}),
\bfpage{2508}--\blpage{2537}
(\byear{2017}).
\doiurl{10.1137/16M1094828}
\end{barticle}
\endbibitem

%%% 36
\bibitem{GeoT'99}
\begin{botherref}
\oauthor{\bsnm{Eberly}, \binits{D.}}:
Least Squares Fitting of Data by Linear or Quadratic Structures,
Redmond WA 98052
(Created: July 15, 1999; Last Modified: September 7, 2021)
\end{botherref}
\endbibitem

%%% 37
\bibitem{snapnets}
\begin{botherref}
\oauthor{\bsnm{Leskovec}, \binits{J.}},
\oauthor{\bsnm{Krevl}, \binits{A.}}:
{SNAP datasets}: {Stanford} large network dataset collection
(2014)
\end{botherref}
\endbibitem

\end{thebibliography}



\end{document}